\documentclass[11pt]{article}

\usepackage{color}
\usepackage{latexsym}
\usepackage{amssymb}
\usepackage{graphicx}
\usepackage{amsmath, amsfonts,amssymb,theorem,euscript,array,enumerate,amsfonts,mathrsfs}

\newtheorem{Theorem}{Theorem}[part]

\newtheorem{Proposition}{Proposition}[part]

\newtheorem{Lemma}{Lemma}[part]

\newtheorem{Remark}{Remark}[part]

\def\esssup_#1{\underset{#1}{\mathrm{ess\,sup\, }}}
\def\essinf_#1{\underset{#1}{\mathrm{ess\,inf\, }}}
\def\argmax_#1{\underset{#1}{\mathrm{arg\,max\, }}}

\def \R{\mathbb{R}}

\def \E{\mathbb{E}}
\def \F{\mathbb{F}}
\def \G{\mathbb{G}}

\def \P{\mathbb{P}}

\def \D{\mathbb{D}}

\def \Ac{{\cal A}}
\def \Bc{{\cal B}}

\def \Dc{{\cal D}}

\def \Fc{{\cal F}}
\def \Gc{{\cal G}}

\def \Pc{{\cal P}}

\def \Oc{{\cal O}}

\def \eps{\varepsilon}

\def \ep{\hbox{ }\hfill$\Box$}

\def\reff#1{{\rm(\ref{#1})}}

\def\beqs{\begin{eqnarray*}}
\def\enqs{\end{eqnarray*}}
\def\beq{\begin{eqnarray}}
\def\enq{\end{eqnarray}}

\begin{document}

\title{Stochastic control under progressive enlargement of filtrations \\ 
and applications to multiple defaults risk management\thanks{I would like to thank N. El Karoui, Y. Jiao and I. Kharroubi for useful comments.}}

\author{  Huy\^en PHAM
             \\\small  Laboratoire de Probabilit\'es et
             \\\small  Mod\`eles Al\'eatoires
             \\\small  CNRS, UMR 7599
             \\\small  Universit\'e Paris 7
             \\\small  e-mail: pham@math.jussieu.fr
             \\\small  and Institut Universitaire de France  
             }


\maketitle

\vspace{20mm}

\begin{abstract}
We formulate and investigate a  general stochastic control problem under a progre\-ssive enlargement 
of filtration. The global information is enlarged from a reference filtration and the knowledge of multiple 
random times together with associated marks when they occur. By working under a density hypothesis on the 
conditional joint distribution of the random times and marks, we prove a  decomposition of the 
original stochastic control problem under the global filtration into  classical stochastic control problems under the reference filtration, 
which are determined in a finite  backward induction.  Our method revisits and extends in particular stochastic control of  diffusion processes with finite number of jumps.   
This  study is  motivated by optimization problems arising in default risk 
management, and we provide applications of our decomposition result for  the indifference pricing of defaultable 
claims, and the optimal investment under bilateral counterparty risk. The solutions are expressed  
in terms of  BSDEs involving only  Brownian filtration,  and remarkably without 
jump terms  coming from the default times and marks  in the global filtration. 
\end{abstract}

\vspace{5mm}


\vspace{3mm}

\noindent {\bf MSC Classification (2000):}  60J75, 93E20, 60H20.

\newpage

\section{Introduction}

The field of stochastic control has known important developments over these last years, inspired especially by various problems in economics and finance arising in risk management, option hedging, optimal investment, portfolio selection or real options valuation.  A vast li\-terature on this topic and its applications has grown with  different approaches  ran\-ging from dynamic programming  method, Hamilton-Jacobi-Bellman Partial Differential Equations (PDEs) and Backward Stochastic Differential Equations (BSDEs) to  convex martingale duality methods. We refer to the 
monographs \cite{fleson06}, \cite{yonzho99}, \cite{okssul07} or \cite{pha08} for recent updates on this subject.  In particular, the theory of BSDEs 
has emerged as a major research topic with original and significant  contributions related to  stochastic control and its financial applications, 
see a recent overview in \cite{elkhammat08}.

On the other hand, the field of enlargement of filtrations is a traditional subject in probability theory  initiated by fundamental works of the  French school in the 80s, see e.g. \cite{jeu80}, \cite{jac85}, \cite{jeuyor85}, and the recent lecture notes \cite{manyor06}.  It knows a renewed interest due to its natural application in credit risk  research  where it appears as a powerful  tool for modelling default events.  For an overview, we refer to the 
books  \cite{bierut02}, \cite{dufsin03},  \cite{sch03} or the lecture notes \cite{biejearut04}. The standard approach of credit  event is based on 
the enlargement of a reference filtration $\F$ (the default-free information structure) by the knowledge of a 
default time when it occurs, leading to the global filtration $\G$, and called progressive enlargement of filtrations. 
Moreover, it assumes that the credit event should arrive by surprise, i.e. it is a totally inacessible 
random time for the reference filtration. Hence, the main approaches consist in modelling the intensity of the random time (usually referred to as 
the reduced-form approach),  or more generally in the modelling of the conditional law of this random time, and referred  to as density hypothesis, see \cite{elkjeajia09}.  The stability of the class of semimartingale, usually called {\bf (H')} hypothesis, and meaning that any $\F$-semimartingale remains a $\G$-semimartingale,  is a fundamental property both in probability and finance where it is closely related to the absence of arbitrage. 
It holds true under the density hypothesis, and the related canonical decomposition in the enlarged filtration can be explicitly expressed,  
as shown in \cite{jeacam09}.  A stronger assumption than {\bf (H')} hypothesis is the so-called immersion property or {\bf (H)} hypothesis, denoting the fact that $\F$-martingales remain $\G$-martingales.

The purpose of this paper is to combine both features of stochastic control and progre\-ssive enlargement of filtrations in view of applications in finance,  in particular for defaults risk management. We formulate and study the general structure for such control pro\-blems by considering 
a progressive enlargement with multiple random times and associated marks. These marks represent for example in credit events  
jump sizes of asset values, which may arrive several times by surprise and cannot be predicted from the past observation of asset  processes.  We work under the density hypothesis on the conditional joint distribution of the random times and marks.  Our new approach consist in decomposing 
the initial control problem in the $\G$-filtration into a finite sequence of control problems formulated in the $\F$-filtration, and which are determined 
recursively.  This is based on an  enlightening re\-presentation of  any $\G$-predictable or optional process that we split into  indexed 
$\F$-predictable or optional processes between each random time.  This point of view allows us to change of regimes in the state process, and to modify  the control set and  the gain functions between random times.  This flexibility in the formulation of the stochastic control problem appears 
also quite useful and re\-levant for financial interpretation.  Our method consist  basically in projecting $\G$-processes into the reference 
$\F$-filtration between two random times,  and features some similarities with  filtering approach. This  contrasts with the standard approach 
in progressive enlargement of filtration  focusing on the  representation of controlled state process  in  the  $\G$-filtration  where the control set has to be  fixed at the initial time. Moreover,  in this global approach, one usually 
assumes that  {\bf (H)} hypothesis holds in order to get a martingale representation in the $\G$-filtration. 
In this case,  the solution is  then characterized  from dynamic programming method in the $\G$-filtration via PDEs with integrodifferential terms  or BSDEs with jumps.  By means of our $\F$-decomposition result under the density hypothesis (and without assuming {\bf (H)} hypothesis), 
we can solve each stochastic control problem by  dynamic programming in the $\F$-filtration, which leads  typically  to PDEs or BSDEs related only to Brownian motion, thus simpler a priori than Integro-PDEs and BSDEs  with jumps.    
Our decomposition method revisits and more importantly extends  stochastic control of diffusion processes with finite number of jumps, and gives 
some new insight  for studying  Integro-PDEs and BSDEs  with jumps.   
We illustrate our methodology with two financial applications in default risk management.  The first one considers  the problem of indifference pricing of defaultable claims, and the second application deals with an optimal 
investment problem under bilateral contagion risk with two nonordered default times. 
The solutions are explicitly expressed in terms of BSDEs  involving only Brownian motion.

The paper is organized as follows.  The next section presents the general framework of progressive enlargement of filtration with 
successive random times and marks.  We  state the decomposition result for a $\G$-predictable and optional process, and as a consequence we derive under the density hypothesis the computation of expectation functionals of $\G$-optional processes in terms of $\F$-expectations. 
In Section 3, we formulate the abstract stochastic control problem in this context and connect it in particular to diffusion processes with jumps. 
Section 4 contains the main $\F$-decomposition result of the initial stochastic control problem.   The case of enlargement of filtration with 
multiple (and not necessarily successive) random times is considered in Section 5, and we show how to derive the results from the case 
of successive random times with auxiliary marks.  Finally,  Section 6 is devoted to 
some applications in risk management, where we  present the results and postpone the detailed proofs and more examples in  a  forthcoming 
paper  \cite{JKP}.

\section{Progressive enlargement of filtration with successive  random times}

\setcounter{equation}{0} \setcounter{Assumption}{0}
\setcounter{Theorem}{0} \setcounter{Proposition}{0}
\setcounter{Corollary}{0} \setcounter{Lemma}{0}
\setcounter{Definition}{0} \setcounter{Remark}{0}

We fix a probability space $(\Omega,\Gc,\P)$, and we start with a reference filtration $\F$ $=$ 
$(\Fc_t)_{t\geq 0}$ satisfying the usual conditions 
($\Fc_0$ contains the null sets of $\P$ and $\F$ is right continuous: $\Fc_t$ $=$ $\Fc_{t^+}$ $:=$ 
$\cap_{s>t}\Fc_s$). We consider a vector of $n$ random times $\tau_1,\ldots,\tau_n$ (i.e. nonnegative $\Gc$-random variables) and  
a vector of $n$ $\Gc$-measurable random variables $\zeta_1,\ldots,\zeta_n$ valued in some Borel subset $E$ of $\R^m$. 
The default information is the knowledge of these  default times $\tau_k$ when they occur, together with the associated marks 
$\zeta_k$. For each $k$ $=$ $1,\ldots,n$, it is defined in mathematical terms as the smallest right-continuous filtration 
$\D^k$ $=$ $(\Dc_t^k)_{t\geq 0}$ such that  $\tau_k$ is a $\D^k$-stopping time, and $\zeta_k$ is $\Dc_{\tau_k}^k$-measurable.  In other words, $\Dc_t^k$ $=$ $\tilde \Dc^k_{t^+}$, where $\tilde \Dc_t^k$ $=$  
$\sigma(\zeta_k 1_{\tau_k\leq s}, s \leq t)$ $\vee$  $\sigma(1_{\tau_k\leq s}, s\leq t)$.  The global market information is then  
defined by the progressive enlargement of filtration $\G$ $=$ $\F$ $\vee$ $\D^1$ $\vee$ $\ldots$ $\vee$ $\D^n$. 
The filtration $\G$ $=$ $(\Gc_t)_{t\geq 0}$  is the smallest filtration containing $\F$, and such that 
for any $k$ $=$ $1,\ldots,n$,  $\tau_k$ is a $\G$-stopping time, and $\zeta_k$ is $\Gc_{\tau_k}$-measurable.  
With respect to  the classical framework of progressive enlargement of filtration with a single random time extensively studied in the literature, 
we consider here multiple random times together with marks. For simplicity of presentation, we first consider the case where the random times 
are ordered, i.e.  $\tau_1\leq\ldots\leq\tau_n$, and so valued  in $\Delta_n$ on $\{\tau_n<\infty\}$, where
\beqs
\Delta_k &=& \Big\{  (\theta_1,\ldots,\theta_k) \in (\R_+)^k: \theta_1\leq \ldots \leq \theta_{k},  \Big\}, \;\;\; k=1,\ldots,n. 
\enqs 
This means actually  that the observations of interest  are the ranked default times (together with the marks). 
We shall indicate in Section 5 how to adapt the results  in the case of multiple random times not necessarily ordered.

\vspace{1mm}

We introduce some  notations  used throughout the paper.  

\vspace{1mm}

\noindent - $\Pc(\F)$ (resp. $\Pc(\G)$) is the $\sigma$-algebra  of $\F$ (resp. $\G$)-predictable measurable subsets on $\R_+\times\Omega$, i.e. 
the $\sigma$-algebra generated by the left-continuous $\F$-adapted (resp. $\G$-adapted) processes.  We also let $\Pc_\F$ (resp. $\Pc_\G$) denote the set of processes that are $\F$-predictable (resp. $\G$-predictable), i.e. $\Pc(\F)$-measurable (resp. $\Pc(\G)$-measurable).

\vspace{1mm}

\noindent - $\Oc(\F)$ (resp. $\Oc(\G)$) is the $\sigma$-algebra  of $\F$ (resp. $\G$)-optional measurable subsets on $\R_+\times\Omega$, i.e. 
the $\sigma$-algebra generated by the right-continuous $\F$-adapted (resp. $\G$-adapted) processes. We also let $\Oc_\F$ (resp. $\Oc_\G$) denote the set of processes that are $\F$-optional (resp. $\G$-optional), i.e. $\Oc(\F)$-measurable (resp. $\Oc(\G)$-measurable).

\vspace{1mm} 
 
\noindent  -  For $k$ $=$ $1,\ldots,n$, we denote by $\Pc^k_{\F}(\Delta_k,E^k)$ (resp. $\Oc^k_{\F}(\Delta_k,E^k)$)  
the set of  of indexed processes $Y^k(.)$ such that the map  
$(t,\omega,\theta_1,\ldots,\theta_k,e_1,\ldots,e_k)$ $\rightarrow$ $Y_t^k(\omega,\theta_1,\ldots,\theta_k,e_1,\ldots,e_k)$ 
is $\Pc(\F)\otimes\Bc(\Delta_k)\otimes\Bc(E^k)$-measurable (resp. $\Oc(\F)\otimes\Bc(\Delta_k)\otimes\Bc(E^k)$-measurable).

\vspace{1mm}

\noindent -  For $\theta$ $=$ $(\theta_1,\ldots,\theta_n)$ $\in$ $\Delta_n$, $e$ $=$ $(e_1,\ldots,e_n)$ $\in$ $E^n$, we denote by
\beqs
\theta^{(k)} \; = \; (\theta_1,\ldots,\theta_k), & &e^{(k)} \; = \; (e_1,\ldots,e_k), \;\;\; k =1,\ldots,n. 
\enqs

The following result provides the key decomposition of  predictable and optional  processes with respect to this progressive enlargement of  filtration. 
This extends a classical result, see e.g. Lemma 4.4 in \cite{jeu80} or Chapter 6 in  \cite{pro06},  
stated for a progressive enlargement of filtration with a single random time.

\begin{Lemma} \label{lemdecG}
Any $\G$-predictable process $Y$ $=$ $(Y_t)_{t\geq 0}$  is represented as
\beq
Y_t &=& Y_t^0 1_{t \leq \tau_1} + \sum_{k=1}^{n-1} Y_t^k(\tau_1,\ldots,\tau_k,\zeta_1,\ldots,\zeta_k) 
1_{\tau_k <  t \leq  \tau_{k+1}}  \nonumber \\
& & \;\;\; + \;  
Y_t^n(\tau_1,\ldots,\tau_n,\zeta_1,\ldots,\zeta_n) 1_{\tau_n <  t}, \;\;\; t \geq 0,  \label{decYpredic}
\enq
where $Y^0$ $\in$ $\Pc_\F$, and $Y^k$ $\in$ $\Pc_\F^k(\Delta_k,E^k)$, for $k$ $=$ $1,\ldots,n$. 
Any $\G$-optional process $Y$ $=$ $(Y_t)_{t\geq 0}$  is represented as
\beq
Y_t &=& Y_t^0 1_{t<\tau_1} + \sum_{k=1}^{n-1} Y_t^k(\tau_1,\ldots,\tau_k,\zeta_1,\ldots,\zeta_k) 
1_{\tau_k\leq t< \tau_{k+1}}  \nonumber \\
& & \;\;\; + \;  
Y_t^n(\tau_1,\ldots,\tau_n,\zeta_1,\ldots,\zeta_n) 1_{\tau_n\leq t}, \;\;\; t \geq 0,  \label{decY}
\enq
where $Y^0$ $\in$ $\Oc_{\F}$, and   $Y^k$ $\in$  $\Oc^k_{\F}(\Delta_k,E^k)$, for $k$ $=$ $1,\ldots,n$.   
\end{Lemma}
{\bf Proof.}  We prove the decomposition result for predictable processes by induction on $n$.  We denote by  $\G^n$ $=$ 
$\F$ $\vee$ $\D^1$ $\vee$ $\ldots$ $\vee$ $\D^n$. 

\noindent {\it Step 1}. Suppose first that $n$ $=$ $1$, so that $\G$ $=$ $\F$ $\vee$ $\D^1$.  Let us  consider  generators of $\Pc(\G)$, which are  processes in the form
\beqs
Y_t &=&  f_s \; g(\zeta_1 1_{\tau_1\leq s}) h(\tau_1\wedge s) 1_{t> s}, \;\;\;\;\; t \geq 0,
\enqs
with $s$ $\geq$ $0$, $f_s$ $\Fc_s$-measurable,  $g$ measurable defined on $E\cup \{0\}$, and $h$ measurable defined on $\R_+$.  By taking 
\beqs
Y_t^0 \; = \;  f_s \; g(0) h(s) 1_{t>s}, & \mbox{ and } & Y_t^1(\theta_1,e) \; = \; f_s g(e 1_{\theta_1\leq s}) h(\theta_1\wedge s) 1_{t> s}, 
\enqs
we see that the decomposition \reff{decYpredic} holds for generators of $\Pc(\G)$.  
We then extend this decomposition for any $\Pc(\G)$-measurable processes, by the monotone class theorem.  

\vspace{1mm}

\noindent {\it Step 2.} Suppose that the result holds for $n$, and consider the case with $n+1$ ranked default times, so that 
$\G$ $=$ $\G^n$ $\vee$ $\D^{n+1}$, $\Dc_t^{n+1}$ $=$ $\tilde \Dc^{n+1}_{t^+}$, where $\tilde \Dc_t^{n+1}$ $=$  
$\sigma(\zeta_{n+1} 1_{\tau_{n+1}\leq s}, s \leq t)$ $\vee$  $\sigma(1_{\tau_{n+1}\leq s}, s\leq t)$. 
 By the same arguments of enlargement of filtration with one  default time as in Step 1, we derive that any 
$\Pc(\G)$-measurable  process $Y$ is represented as
\beq \label{decYinter}
Y_t &=& Y_t^{0,(n)} 1_{t\leq\tau_{n+1}} + Y_t^{1,(n)}(\tau_{n+1},\zeta_{n+1}) 1_{\tau_{n+1}<t},
\enq
where $Y^{0,(n)}$ is $\Pc(\G^n)$-measurable, and $(t,\omega,\theta_{n+1},e_{n+1})$ $\mapsto$ $Y_t^{1,(n)}(\omega,\theta_{n+1},e_{n+1})$ is 
$\Pc(\G^n)\otimes\Bc(\R_+)\otimes\Bc(E)$-measurable.  Now, from the induction hypothesis for $\G^n$, we have
\beqs
Y_t^{0,(n)} &=& Y_t^{0,0,(n)} 1_{t \leq \tau_1} + \sum_{k=1}^{n-1} Y_t^{k,0,(n)} (\tau_1,\ldots,\tau_k,\zeta_1,\ldots,\zeta_k) 
1_{\tau_k <  t \leq  \tau_{k+1}}  \nonumber \\
& & \;\;\; + \;  
Y_t^{n,0,(n)}(\tau_1,\ldots,\tau_n,\zeta_1,\ldots,\zeta_n) 1_{\tau_n <  t}, \;\;\; t \geq 0, 
\enqs
where $Y^{0,0,(n)}$ $\in$ $\Pc_\F$, and $Y^{k,0,(n)}$ $\in$ $\Pc_\F^k(\Delta_k,E^k)$, for $k$ $=$ $1,\ldots,n$. Similarly, we have 
\beqs
Y_t^{1,(n)}(\theta_{n+1},e_{n+1}) &=& Y_t^{0,1,(n)}(\theta_{n+1},e_{n+1}) 1_{t \leq \tau_1} \\
& & \;\;\; + \;  \sum_{k=1}^{n-1} Y_t^{k,1,(n)} (\tau_1,\ldots,\tau_k,\zeta_1,\ldots,\zeta_k,\theta_{n+1},e_{n+1}) 
1_{\tau_k <  t \leq  \tau_{k+1}}  \nonumber \\
& & \;\;\; + \;  
Y_t^{n,1,(n)}(\tau_1,\ldots,\tau_n,\zeta_1,\ldots,\zeta_n,\theta_{n+1},e_{n+1}) 1_{\tau_n <  t}, \;\;\; t \geq 0, 
\enqs
where $Y^{0,1,(n)}$ $\in$ $\Pc_\F^1(\R_+,E)$,  $Y^{k,1,(n)}$ $\in$ $\Pc_\F^{k+1}(\Delta_k\times\R_+,E^{k+1})$, $k$ $=$ $1,\ldots,n$. 
Finally, plugging these two decompositions with respect to $\Pc(\G^n)$ into relation \reff{decYinter}, and recalling that 
$\tau_1\leq\ldots\leq\tau_n\leq\tau_{n+1}$, we get the required decomposition at level $n+1$ for $\G$: 
\beqs
Y_t &=&  Y_t^{0,0,(n)} 1_{t \leq \tau_1} + \sum_{k=1}^{n} Y_t^{k,0,(n)} (\tau_1,\ldots,\tau_k,\zeta_1,\ldots,\zeta_k) 
1_{\tau_k <  t \leq  \tau_{k+1}} \\
& & \;\;\; + \;  
Y_t^{n+1}(\tau_1,\ldots,\tau_{n+1},\zeta_1,\ldots,\zeta_{n+1}) 1_{\tau_{n+1} <  t}, \;\;\; t \geq 0, 
\enqs
where we notice that the indexed process $Y^{n+1}$ defined by  $Y^{n+1}(\theta_1,\ldots,\theta_{n+1},e_1,\ldots,e_{n+1})$ $:=$ 
$Y^{n,1,(n)}(\theta_1,\ldots,\theta_n,e_1,\ldots,e_n,\theta_{n+1},e_{n+1})$,  lies in $\Pc^{n+1}_\F(\Delta_{n+1},E^{n+1})$.

The decomposition result for $\G$-optional processes is proved similarly by induction and considering generators of $\Oc(\G^1)$, 
which are  processes in the form
\beqs
Y_t &=&  f_s \; g(\zeta_1 1_{\tau_1\leq s}) h(\tau_1\wedge s) 1_{t \geq s}, \;\;\;\;\; t \geq 0,
\enqs
with $s$ $\geq$ $0$, $f_s$ $\Fc_s$-measurable,  $g$ measurable defined on $E\cup \{0\}$, and $h$ measurable defined on $\R_+$. 
\ep

\vspace{2mm}

In view of the decomposition \reff{decYpredic} or  \reff{decY},  we can then identify any  $Y$  $\in$ $\Pc_{\G}$ (resp. $\Oc_\G$) with an $n+1$-tuple 
$(Y^0,\ldots,Y^n)$ $\in$ $\Pc_{\F}\times\ldots\times\Pc^n_{\F}(\Delta_n,E^n)$ (resp.  $\Oc_{\F}\times\ldots\times\Oc^n_{\F}(\Delta_n,E^n)$).

\vspace{2mm}

We now  require a  density hypothesis  on the random times and their associated jumps by assuming that  for any $t$, the conditional 
distribution  of $(\tau_1,\ldots,\tau_n,\zeta_1,\ldots,\zeta_n)$ given  $\Fc_t$ is absolutely continuous with respect to a positive 
measure $\lambda(d\theta) \eta(de)$ on $\Bc(\Delta_n)\otimes\Bc(E^n)$, with $\lambda$  the Lebesgue measure $\lambda(d\theta)$ $=$ 
$d\theta_1\ldots d\theta_n$, and $\eta$ a product measure  $\eta(de)$ $=$ $\eta_1(de_1)\ldots\eta_1(de_n)$ on $\Bc(E)\otimes\ldots\otimes\Bc(E)$. 
More precisely, we assume that  there exists $\gamma$ $\in$ $\Oc^n_{\F}(\Delta_n,E^n)$  such that
\beqs
{\bf (DH)} \hspace{7mm} 
\P\big[ (\tau_1,\ldots,\tau_n,\zeta_1,\ldots,\zeta_n) \in d\theta  de | \Fc_t \big] 
&=& \gamma_t(\theta_1,\ldots,\theta_n,e_1,\ldots,e_n)  \\
& & d\theta_1\ldots d\theta_n \; \eta_1(de_1)\ldots\eta_1(de_n), \;\;\; a.s. 
\enqs

\begin{Remark}
{\rm  In the particular case where $\gamma$ is in the form $\gamma_t(\theta,e)$ $=$ $\varphi_t(\theta)\psi_t(e)$, 
the condition {\bf (DH)}  means  that the random times $(\tau_1,\ldots,\tau_n)$ and the jump sizes $(\zeta_1,\ldots,\zeta_n)$ are 
independent given $\Fc_t$,  for all $t$ $\geq$ $0$, and 
\beqs  
\P\Big[ (\tau_1,\ldots,\tau_n) \in d\theta  | \Fc_t \big] \; = \;  \varphi_t(\theta) \lambda(d\theta), & & 
\P\Big[ (\zeta_1,\ldots,\zeta_n) \in de  | \Fc_t \big] \; = \;  \psi_t(e) \eta(de), \;\;\; a.s. 
\enqs
This condition extends the usual  density hypothesis for a random time in the theory of initial or progressive enlargement of filtration, see  
\cite{jac85} or  \cite{jeacam09}.  An important  result in the theory of enlargement of filtration under the density hypothesis is the semimartingale invariance property, also called {\bf (H')} hypothesis, i.e.  any $\F$-semimartingale remains a $\G$-semimartingale. This result is related in finance to no-arbitrage conditions, and is thus also a desirable property from a economical viewpoint. Random times satisfying the density hypothesis are very well suitable for the analysis of credit risk   events, as shown  recently in  \cite{elkjeajia09}.  We also refer to this paper for a  discussion 
on the relation between the density hypothesis and the reduced-form (or intensity) approach in credit risk modelling. 
}
\end{Remark}

\vspace{2mm}

In the sequel, it is useful  to introduce the following notations.  We denote by $\gamma^0$  the $\F$-optional process defined by 
\beq
\gamma_t^0 &=& \P\big[ \tau_1 > t \big| \Fc_t \big]  \label{surgamma0} \\
&=& \int_{E^n}  \int_t^\infty \int_{\theta_1}^\infty \ldots \int_{\theta_{n-1}}^\infty 
\gamma_t(\theta_1,\ldots,\theta_n,e_1,\ldots,e_n) d\theta_n\ldots d\theta_1 \eta_1(de_1)\ldots\eta_n(de_n), \nonumber
\enq
and  we denote by  $\gamma^k$,  $k$ $=$ $1,\ldots,n-1$,  the indexed process in $\Oc_{\F}^k(\Delta_k,E^k)$  defined by 
\beqs
& & \gamma_t^k(\theta_1,\ldots,\theta_k,e_1,\ldots,e_k) \\
&=& \int_{E^{n-k}} \int_t^\infty  \int_{\theta_{k+1}}^\infty \ldots \int_{\theta_{n-1}}^\infty \gamma_t(\theta_1,\ldots,\theta_n,e_1,\ldots,e_n) 
d\theta_n\ldots d\theta_{k+1} \eta_1(de_{k+1})\ldots\eta_1(de_n),
\enqs
so that for $k$ $=$ $1,\ldots,n-1$, 
\beq \label{surgammak}
\P\big[ \tau_{k+1} > t  \big| \Fc_t\big] &=&  \int_{E^k} \int_{\Delta_k} \gamma_t^k(\theta_1,\ldots,\theta_k,e_1,\ldots,e_k) 
d \theta_1 \ldots d\theta_k \eta_1(de_1)\ldots\eta_1(de_k). 
\enq
Notice that the family of measurable maps $\gamma_k$, $k$ $=$ $0,\ldots,n$ can be also written in backward induction by
\beqs
\gamma_t^k(\theta_1,\ldots,\theta_k,e_1,\ldots,e_k) &=&  
\int_E \int_t^\infty \gamma_t^{k+1} (\theta_1,\ldots,\theta_{k+1},e_1,\ldots,e_{k+1}) d\theta_{k+1} \eta_{1}(de_{k+1}), 
\enqs
for $k$ $=$ $0,\ldots,n-1$, starting from $\gamma^n$ $=$ $\gamma$.  In view of \reff{surgamma0}-\reff{surgammak}, the process $\gamma_k$ may be interpreted as the survival density process of $\tau_{k+1}$, $k$ $=$ $0,\ldots,n-1$.


\vspace{2mm}

The next result  provides the computation for the optional projection of a $\Oc(\G)$-measurable process on  the reference filtration $\F$.

\begin{Lemma} \label{lemopt}
Let $Y$ $=$ $(Y^0,\ldots,Y^n)$ be a nonnegative (or bounded)  $\G$-optional process. Then for any $t$ $\geq$ $0$, we have
\beqs
\hat Y_t^{\F} &:=& \E\big[ Y_t \big| \Fc_t \big] \\
&=& Y_t^0 \gamma_t^0 + \sum_{k=1}^{n} \int_{E^{k}} \int_0^t\ldots\int_{\theta_{k-1}}^t 
Y_t^{k}(\theta_1,\ldots,\theta_{k},e_1,\ldots,e_{k}) \\
& &   \hspace{3cm}  \gamma_t^{k}(\theta_1,\ldots,\theta_k,e_1,\ldots,e_k) d\theta_k\ldots d\theta_1 \eta_1(de_1)\ldots\eta_1(de_k),
\enqs
where we used the convention that $\theta_{k-1}$ $=$ $0$ for $k$ $=$ $1$ in the above integral.  Equivalently, we have the backward induction 
formula for $\hat Y_t^{\F}$ $=$ $\hat Y_t^{0,\F}$, where the $\hat Y_t^{k,\F}$ are given  for any $t$ $\geq$ $0$, by
\beqs
\hat Y_t^{n,\F}(\theta,e) &=& 
Y_t^n(\theta,e)\gamma_t(\theta,e) \\
\hat Y_t^{k,\F}(\theta^{(k)},e^{(k)}) &=&  
Y_t^k(\theta^{(k)},e^{(k)}) \gamma_t^k(\theta^{(k)},e^{(k)}) \\
& & \; + \; \int_{\theta_k}^t  \int_E \hat Y_t^{k+1,\F}(\theta^{(k)},\theta_{k+1},e^{(k)},e_{k+1}) \eta_1(de_{k+1}) d \theta_{k+1}, 
\enqs
for  $\theta$ $=$ $(\theta_1,\ldots,\theta_n)$ $\in$ $\Delta_n\cap [0,t]^n$, $e$ $=$ $(e_1,\ldots,e_n)$ $\in$ $E^n$. 
\end{Lemma}
{\bf Proof.}
Let $Y$ $=$ $(Y^0,\ldots,Y^n)$ be a nonnegative (or bounded)  $\G$-optional process, decomposed as in \reff{decY} so that: 
\beq \label{Yprojinter}
\E\big[ Y_t \big| \Fc_t \big]  &=& \E\Big[ Y_t^0 1_{t<\tau_1} \Big| \Fc_t\Big]  
+ \sum_{k=1}^n \E\Big[ Y_t^k(\tau_1,\ldots,\tau_k,\zeta_1,\ldots,\zeta_k)  1_{\tau_k\leq t< \tau_{k+1}} \Big| \Fc_t \Big], 
\enq
with the convention that $\tau_{n+1}$ $=$ $\infty$. 
Now, for any $k$ $=$ 
$1,\ldots,n$, we have under the density hypothesis {\bf (DH)}
\beqs
& & \E\Big[ Y_t^k(\tau_1,\ldots,\tau_k,\zeta_1,\ldots,\zeta_k)  1_{\tau_k\leq t< \tau_{k+1}} \Big| \Fc_t \Big]  \\
&=&  \int_{\Delta_n\times E^n} Y_t^k(\theta_1,\ldots,\theta_k,e_1,\ldots,e_k) 1_{\theta_k \leq t<\theta_{k+1}} 
\gamma_t(\theta_1,\ldots,\theta_n,e_1,\ldots,e_n) \lambda(d\theta) \eta(de) \\
&=& \int_{E^{k}} \int_0^t\ldots\int_{\theta_{k-1}}^t  Y_t^{k}(\theta_1,\ldots,\theta_{k},e_1,\ldots,e_{k}) 
 \gamma_t^{k}(\theta_1,\ldots,\theta_k,e_1,\ldots,e_k) d\theta_k\ldots d\theta_1 \\
 & & \hspace{11cm} \eta_1(de_1)\ldots\eta_1(de_k),
\enqs
where the second inequality follows from Fubini's theorem and the definition of $\gamma^k$.  We also have
\beqs
\E\Big[ Y_t^0 1_{t<\tau_1} \Big| \Fc_t\Big] &=& Y_t^0 \P[ \tau_1 > t | \Fc_t] \; = \; Y_t^0 \gamma_t^0. 
\enqs
We then get the required result by plugging these two last relations into \reff{Yprojinter}.  Finally, the backward formula for the $\F$-optional projection 
of $Y$ is  obtained by a straightforward induction. 
\ep

\vspace{2mm}

As a consequence of the above backward induction formula for the optional projection, we derive a bakward  formula for the computation of 
expectation functionals of  $\G$-optional processes, which involves only $\F$-expectations.

\begin{Proposition} \label{corollexp}
Let $Y$ $=$ $(Y^0,\ldots,Y^n)$  and $Z$ $=$ $(Z^0,\ldots,Z^n)$ be two nonnegative (or bounded)  $\G$-optional processes,  
and fix $T$ $\in$ $(0,\infty)$.

\noindent The expectation $\E[\int_0^T Y_t dt + Z_T]$ can be computed in a backward induction as 
\beqs
\E\Big[\int_0^T Y_t dt + Z_T\Big] &=&  J_0
\enqs
where the  $J_k$, $k$ $=$ $0,\ldots,n$ are given by
\beqs
J_n(\theta,e) &=& \E \Big[ \int_{\theta_n}^T Y_t^n \gamma_t(\theta,e) dt + Z_T^n  \gamma_T(\theta,e) \Big| \Fc_{\theta_n} \Big] \\
J_k(\theta^{(k)},e^{(k)}) &=&   \E \Big[  \int_{\theta_k}^T Y_t^k  \gamma_t^k(\theta^{(k)},e^{(k)}) dt + 
Z_T^k  \gamma_T^k(\theta^{(k)},e^{(k)}) \\
& & \;\;\; + \;  \int_{\theta_k}^T  \int_E  J_{k+1}(\theta^{(k)},\theta_{k+1},e^{(k)},e_{k+1}) \eta_1(de_{k+1}) d \theta_{k+1} \Big| \Fc_{\theta_k} \Big], 
\enqs
for  $\theta$ $=$ $(\theta_1,\ldots,\theta_n)$ $\in$ $\Delta_n\cap [0,T]^n$, $e$ $=$ $(e_1,\ldots,e_n)$ $\in$ $E^n$, with the convention $\theta_0$ 
$=$ $0$. 
\end{Proposition}
{\bf Proof.}  For any  $\theta$ $=$ $(\theta_1,\ldots,\theta_n)$ $\in$ $\Delta_n\cap [0,T]^n$, $e$ $=$ $(e_1,\ldots,e_n)$ $\in$ 
$E^n$, let us define 
\beqs
J_k(\theta^{(k)},e^{(k)}) &=& \E\Big[\int_{\theta_k}^T \hat Y_t^{k,\F}(\theta^{(k)},e^{(k)}) dt + \hat Z_T^{k,\F} (\theta^{(k)},e^{(k)}) \big|\Fc_{\theta_k}\Big],
\enqs 
where the $\hat Y^{k,\F}$ and   $\hat Z^{k,\F}$ are defined in  Lemma \ref{lemopt}, associated respectively to $Y$ and $Z$.  
Then $J_0$ $=$ $\E[\int_0^T Y_t dt + Z_T]$, and 
 we see from the backward induction for $\hat Y^{k,\F}$ and   $\hat Z^{k,\F}$  that the $J_k$, $k$ $=$ $0,\ldots,n$,  satisfy 
\beqs
J_n(\theta,e) &=& \E \Big[ \int_{\theta_n}^T Y_t^n \gamma_t(\theta,e) dt + Z_T^n(\theta,e) \gamma_T(\theta,e) \Big| \Fc_{\theta_n} \Big] \\
J_k(\theta^{(k)},e^{(k)}) &=&   \E \Big[  \int_{\theta_k}^T Y_t^k  \gamma_t^k(\theta^{(k)},e^{(k)}) dt  + 
Z_T^k  \gamma_T^k(\theta^{(k)},e^{(k)}) \\
& & \;\;\; + \; \int_{\theta_k}^T \int_{\theta_k}^t  \int_E \hat Y_t^{k+1,\F}(\theta^{(k)},\theta_{k+1},e^{(k)},e_{k+1}) \eta_1(de_{k+1}) d \theta_{k+1} dt \\
& & \;\;\; + \;  \int_{\theta_k}^T  \int_E \hat Z_T^{k+1,\F}(\theta^{(k)},\theta_{k+1},e^{(k)},e_{k+1}) \eta_1(de_{k+1}) d \theta_{k+1}
\Big| \Fc_{\theta_k} \Big] \\
&=&  \E \Big[  \int_{\theta_k}^T Y_t^k  \gamma_t^k(\theta^{(k)},e^{(k)}) dt  
+ Z_T^k  \gamma_T^k(\theta^{(k)},e^{(k)})\\
& & \;\;\; + \; \int_{\theta_k}^T \int_E \int_{\theta_{k+1}}^T   \hat Y_t^{k+1,\F}(\theta^{(k)},\theta_{k+1},e^{(k)},e_{k+1}) dt \; \eta_1(de_{k+1}) 
d \theta_{k+1}  \\
& & \;\;\; + \;  \int_{\theta_k}^T  \int_E \hat Z_T^{k+1,\F}(\theta^{(k)},\theta_{k+1},e^{(k)},e_{k+1}) \eta_1(de_{k+1}) d \theta_{k+1}
\Big| \Fc_{\theta_k} \Big] \\
&=& \E \Big[  \int_{\theta_k}^T Y_t^k  \gamma_t^k(\theta^{(k)},e^{(k)}) dt + Z_T^k  \gamma_T^k(\theta^{(k)},e^{(k)})\\
& & \;\;\; + \; \int_{\theta_k}^T \int_E  J_{k+1} (\theta^{(k)},\theta_{k+1},e^{(k)},e_{k+1})   \; \eta_1(de_{k+1}) d \theta_{k+1} 
\Big| \Fc_{\theta_k} \Big],
\enqs
where we used Fubini's theorem in the second equality for $J_k$, and the law of iterated condional expectations for the last equality. This proves the required induction formula for $J_k$, $k$ $=$ $0,\ldots,n$. 
\ep

\section{Abstract stochastic  control problem}

\setcounter{equation}{0} \setcounter{Assumption}{0}
\setcounter{Theorem}{0} \setcounter{Proposition}{0}
\setcounter{Corollary}{0} \setcounter{Lemma}{0}
\setcounter{Definition}{0} \setcounter{Remark}{0}

In this section, we formulate  the  general stochastic control problem in the context of progressively enlargement of filtration with successive 
random times and marks.

\subsection{Controls and state process}

A control is a $\G$-predictable process $\alpha$ $=$ $(\alpha^0,\ldots,\alpha^n)$ $\in$ $\Pc_{\F}\times \ldots\times \Pc^n_{\F}(\Delta_n,E^n)$, 
where  the $\alpha^k$, $k=0,\ldots,n$,   are valued  in some given Borel set $A_k$  of an Euclidian space.  
We denote by $\Pc_\F(A_0)$ (resp. $\Pc_\F^k(\Delta_k,E^k;A_k)$, $k$ $=$ $1,\ldots,n$),  the set of elements in $\Pc_\F$ 
(resp. $\Pc_\F^k(\Delta,E^k)$, $k$ $=$ $1,\ldots,n$)  valued in $A_0$ (resp. $A_k$, $k$ $=$ $1,\ldots,n$).  We set $A$ $=$ 
$A_0\times\ldots\times A_n$, and denote by  $\Ac_\G$ the set of {\it admissible controls} as the product $\Ac_\F^0\times\ldots\times\Ac_\F^n$, where 
$\Ac_\F^0$ (resp. $\Ac_\F^k$, $k$ $=$ $1,\ldots,n$)  is some separable metric space of $\Pc_\F(A_0)$ (resp. $\Pc_\F^k(\Delta_k,E^k;A_k)$, $k$ $=$ $1,\ldots,n$).  The separability condition is required for measurability selection issue.

\vspace{2mm}

The description of the controlled state process is  formulated as follows: 

\noindent $\bullet$ {\it Controlled state process between default times:} we are given a collection of measurable mappings: 
\beq
(x,\alpha^0) \in \R^d\times\Ac^0_{\F} & \longmapsto &  X^{0,x,\alpha^0} \; \in \;  \Oc_{\F} \label{X0con} \\
(x,\alpha^k) \in \R^d\times\Ac^k_{\F}  & \longmapsto &  X^{k,x,\alpha^k} \; \in \;  \Oc^k_{\F}(\Delta_k,E^k), \;\;\; k =1,\ldots,n, \label{Xkcon}
\enq
such that  
we have the initial data: 
\beqs
X_0^{0,x,\alpha^0} &=& x, \;\; \forall x \in \R^d,  \\
X_{\theta_k}^{k,\xi,\alpha^k}(\theta_1,\ldots,\theta_k,e_1,\ldots,e_k) &=& \xi, \;\;\; \forall \xi  \;  \Fc_{\theta_k}-\mbox{measurable}, \;  k =1,\ldots,n. 
\enqs

\noindent $\bullet$ {\it  Jumps of the controlled state process:}  we are  given a collection of maps $\Gamma^k$ on $\R_+\times\Omega\times\R^d\times A_{k-1}\times E$, for $k$ $=$ $1,\ldots,n$,  
such that
\beqs
(t,\omega,x,a,e) & \mapsto & \Gamma^k_t(\omega,x,a,e) \; \mbox{ is } 
\Pc(\F)\otimes\Bc(\R^d)\otimes\Bc(A_{k-1})\otimes\Bc(E)-\mbox{measurable}. 
\enqs

\noindent $\bullet$ {\it Global controlled state process:}  the controlled state process is then given by the mapping 
\beqs
(x,\alpha=(\alpha^0,\ldots,\alpha^n)) \in \R^d\times\Ac_{\G} & \longmapsto & X^{x,\alpha} \in \Oc_{\G},
\enqs
where $X^{x,\alpha}$ is the process equal to
\beq
X_t^{x,\alpha} &=&  \bar X_t^0 1_{t<\tau_1} + \sum_{k=1}^{n-1} \bar X_t^k(\tau_1,\ldots,\tau_k,\zeta_1,\ldots,\zeta_k) 
1_{\tau_k\leq t< \tau_{k+1}}  \nonumber \\
& & \;\;\; + \;   \bar X_t^n(\tau_1,\ldots,\tau_n,\zeta_1,\ldots,\zeta_n) 1_{\tau_n\leq t}, \;\;\; t \geq 0, \label{dynXcon}
\enq
with $(\bar X^0,\ldots,\bar X^n)$ $\in$ $\Oc_{\F}\times\ldots\times\Oc^n_{\F}(\Delta_n,E^n)$  given by
\beqs
\bar X^0 &=& X^{0,x,\alpha^0} \\
\bar X^k(\theta_1,\ldots,\theta_k,e_1,\ldots,e_k) &=& 
X^{k,\Gamma^k_{\theta_k}(\bar X^{k-1}_{\theta_k},\alpha^{k-1}_{\theta_k},e_k),\alpha^k}(\theta_1,\ldots,\theta_k,e_1,\ldots,e_k),
\enqs
for $k$ $=$ $1,\ldots,n$.

The interpretation is the following.  Between the time interval $\tau_k$ $=$ $\theta_k$ and $\tau_{k+1}$ $=$ $\theta_{k+1}$,  
 $k$ $=$ $0,\ldots,n-1$ (with the convention $\theta_0$ $=$ $0$),  the state process $X$ $=$ $\bar X^k$  is  controlled  
 by  $\alpha_k$, which is based  on the basic information $\F$, and  the knowledge of the past jump times and marks 
 $(\theta_1,\ldots,\theta_k,e_1,\ldots,e_k)$.   
Then, at time $\theta_{k+1}$, there is a jump on the state process  determined by the map $\Gamma^{k+1}$, which depends on the 
current state value, control and information, but also on a  ``nonpredictable"  mark $\zeta_{k+1}$ $=$ $e_{k+1}$ at time $\theta_{k+1}$: 
\beqs
X_{\tau_{k+1}} &=& \Gamma_{\tau_{k+1}}^{k+1}(X_{\tau_{k+1}^-},\alpha_{\tau_{k+1}}^k,\zeta_{k+1}). 
\enqs

\subsection{Typical controlled state process}

In typical applications, the dynamics of $X^0$ $=$ $X^{0,x,\alpha^0}$, $X^k$ $=$ $X^{k,x,\alpha^k}$, $k$ $=$ $1,\ldots,n$, 
are governed by diffusion processes: 
\beq
d X_t^0 &=&  b^0_t(X_t^0,\alpha_t^0) dt + \sigma_t^0(X_t^0,\alpha_t^0) dW_t, \;\;\; t \geq 0  \label{dynX0}  \\
d X_t^k  &=& b^k_t(X_t^k,\alpha^k_t,\theta_1,\ldots,\theta_k,e_1,\ldots,e_k) dt  \label{dynXk} \\
& & \;\;\; + \;  \sigma^k_t(X_t^k,\alpha_t^k,\theta_1,\ldots,\theta_k,e_1,\ldots,e_k) dW_t, \;\;\;   t \geq \theta_k, \nonumber 
\enq
Here, $W$ is a standard $m$-dimensional $(\P,\F)$-Brownian motion, and  $(t,\omega,x,a)$ $\rightarrow$ $b_t^0(\omega,x,a)$,  
$\sigma_t^0(\omega,x,a)$ are $\Pc(\F)\otimes\Bc(\R^d)\otimes\Bc(A_0)$-measurable maps valued respectively in $\R^d$ and $\R^{d\times m}$,   
for $k$ $=$ $1,\ldots,n$,  the maps $(t,\omega,x,a,\theta_1,\ldots,\theta_k,e_1,\ldots,e_k)$ $\rightarrow$ 
$b_t^k(\omega,x,a,\theta_1,\ldots,\theta_k,e_1,\ldots,e_k)$, $\sigma_t^k(\omega,x,u,\theta_1,\ldots,\theta_k,e_1,\ldots,e_k)$ are 
$\Pc(\F)\otimes\Bc(\R^d)\otimes\Bc(A_k)\otimes\Bc(\Delta_k)\otimes\Bc(E^k)$-measurable valued respectively in $\R^d$ and $\R^{d\times m}$.    
To alleviate notations, we omitted in  \reff{dynXk}   the dependence of $X^k$, $\alpha^k$ in 
$(\theta_1,\ldots,\theta_k,e_1,\ldots,e_k)$.  We make the linear growth and Lipschitz assumptions on  the functions $x$ $\rightarrow$ 
$b_t^k(x,.)$, $\sigma^k(x,.)$, $k$ $=$ $0,\ldots,n$, in order to ensure for all $(\theta_1,\ldots,\theta_k,e_1,\ldots,e_k)$ $\in$ $\Delta_k\times E^k$, 
the existence and uniqueness of a solution $X^k(\theta_1,\ldots,\theta_k,e_1,\ldots,e_k)$ to the sde \reff{dynX0}, \reff{dynXk}, given the controls and the initial conditions, and this indexed process $X^k$ lies in $\Oc_\F^k(\Delta_k,E^k)$.   The dependence of  the coefficients 
$b^k$, $\sigma^k$ on  the past jump times $\theta_1,\ldots,\theta_k$, and marks $e_1,\ldots,e_k$,  
corresponds to  change of regimes after each jump time, and may be interpreted  in finance as  rating  upgrades or downgrades.   Also, 
typical choice for the set of admissible controls $\Ac_{\F}^k$ is subset  of  indexed $\F$-predictable processes in $L^p$, $p$ $\in$ $[1,\infty)$, 
and the separability of $\Ac_{\F}^k$ follows from the separability of $L^p$, see the discussion in \cite{sontou02}. 

\vspace{3mm}

\noindent {\bf Connection with controlled jump-diffusion processes.}
 
\noindent {\rm Consider the particular  case where  the sets  of controls $A_k$ are identical, equal to $A$, and  let us define the mappings 
$b$ and $\sigma$ on $\R_+\times\Omega\times\R^d\times A$ by: 
\beqs
b_t(x,a) &=& b_t^0(x,a) 1_{t\leq \tau_1} + \sum_{k=1}^{n-1} b_t^k(x,a,\tau_1,\ldots,\tau_k,\zeta_1,\ldots,\zeta_k) 
1_{\tau_k <  t \leq \tau_{k+1}} \\
& & \;\;\; + \;  b_t^n(x,a,\tau_1,\ldots,\tau_n,\zeta_1,\ldots,\zeta_n)  1_{t>\tau_{n}}, \\
\sigma_t(x,a) &=& \sigma_t^0(x,a) 1_{t\leq \tau_1} + \sum_{k=1}^{n-1} \sigma_t^k(x,a,\tau_1,\ldots,\tau_k,\zeta_1,\ldots,\zeta_k) 
1_{\tau_k <  t \leq \tau_{k+1}} \\
& & \;\;\; + \;  \sigma_t^n(x,a,\tau_1,\ldots,\tau_n,\zeta_1,\ldots,\zeta_n)  1_{t > \tau_{n}},
\enqs
and notice that  the maps  $(t,\omega,x,a)$ $\rightarrow$ $b_t(\omega,x,a)$, $\sigma_t(\omega,x,a)$ are 
$\Pc(\G)\otimes\Bc(\R^d)\otimes\Bc(A)$-measurable.  Denote also by $\delta$ the mapping on 
$\R_+\times\Omega\times\R^d\times A\times E$:
\beqs
\delta_t(x,a,e) &=& \sum_{k=1}^{n-1}\Big( \Gamma^k_t(x,a,e) - x\Big) 1_{\tau_k < t \leq  \tau_{k+1}} + 
\Big( \Gamma^n_t(x,a,e) - x\Big)  1_{t>\tau_n},
\enqs
which is  $\Pc(\G)\otimes\Bc(\R^d)\otimes\Bc(A)\otimes\Bc(E)$-measurable. Let us denote by $\mu(dt,de)$  the integer-valued random measure 
associated to the times $\tau_k$ and the marks $\zeta_k$, $k$ $=$ $1,\ldots,n$, which is then given by
\beqs
\mu([0,t]\times B) &=& \sum_{k\geq 1} 1_{\tau_n\leq t} 1_B(\zeta_k), \;\;\; \forall t \geq 0, \; B \in \Bc(E). 
\enqs
The progressive enlarged filtration $\G$ can then be written also as: $\G$ $=$ $\F$ $\vee$ $\F^\mu$ where $\F^\mu$ 
is the right-continuous filtration generated by the integer-valued random measure $\mu$. 
Now, since  the semimartingale property is preserved under  the density hypothesis for this progressive enlargement of filtration, 
(see  \cite{jeacam09}), the process  $W$ remains a semimartingale under $(\P,\G)$ (with a canonical decomposition, which can be explicitly expressed in terms of the density).  Then, we can write the dynamics 
of the state process $X$ $=$ $X^{x,\alpha}$ in \reff{dynXcon} as a controlled jump-diffusion process under  $(\P,\G)$: 
\beqs
dX_t &=& b_t(X_t,\alpha_t) dt + \sigma_t(X_t,\alpha_t) dW_t  + \int_E \delta_t(X_{t^-},\alpha_t,e) \mu(dt,de). 
\enqs
However, notice that in the above $\G$-formulation, the process $W$ is not in  general a Brownian motion under $(\P,\G)$, unless the so-called 
{\bf (H)} immersion property is satisfied, i.e. the martingale property is preserved  from $\F$ to $\G$, which corresponds to the particular case where the density satisfies: $\gamma_t(\theta,e)$ $=$ $\gamma_\theta(\theta,e)$ for $t$ $\geq$ $\theta$.

In the classical formulation by controlled jump-diffusion processes,  one has to fix  a control set $A$, which is invariant during the time  horizon. 
Here, the more general formulation \reff{dynXcon}  allows us to consider different control sets $A_k$ between two default times, and this may be relevant in practical applications. Moreover, we have a suitable decomposition of the coefficients and controlled state process between random times, which provides a natural interpretation in economics and finance.

\subsection{Stochastic control problem}

 In the general framework  for the controlled process in \reff{dynXcon},  let us formulate the objective function for the stochastic control problem on a finite horizon $T$. The terminal gain function is given by  a nonnegative map $G_T$ on $\Omega\times\R^d$ such that $(\omega,x)$ 
 $\mapsto$ $G_T(\omega,x)$ is 
 $\Gc_T\otimes\Bc(\R^d)$-measurable, and which may be represented as
 \beqs
 G_T(x) &=& G_T^0(x) 1_{T<\tau_1} + \sum_{k=1}^{n-1} G_T^k(x,\tau_1,\ldots,\tau_k,\zeta_1,\ldots,\zeta_k) 
 1_{\tau_k\leq T<\tau_{k+1}} \\
 & & \;\;\; + \; G_T^n(x,\tau_1,\ldots,\tau_n,\zeta_1,\ldots,\zeta_n) 1_{\tau_n\leq T},
 \enqs
 where $G_T^0$ is $\Fc_T\otimes\Bc(\R^d)$-measurable, and 
 $G_T^k$ is $\Fc_T\otimes\Bc(\R^d)\otimes\Bc(\Delta_k)\otimes\Bc(E^k)$-measurable, for  $k$ $=$ $1,\ldots,n$.  
 The running gain function is given by a nonnegative map $f$ on $\Omega\times\R^d\times A$ such that 
 $(t,\omega,x,a)$ $\mapsto$ $f_t(\omega,x,a)$ is $\Oc(\G)\otimes\Bc(\R^d)\otimes\Bc(A)$-measurable, and which may be decomposed as 
 \beqs
 f_t(x,a) &=& f_t^0(x,a_0) 1_{t < \tau_1} +  \sum_{k=1}^{n-1} f_t^k(x,a_k,\tau_1,\ldots,\tau_k,\zeta_1,\ldots,\zeta_k) 
 1_{\tau_k\leq t <\tau_{k+1}} \\
 & &   \;\;\; + \; f_t^n(x,a_n,\tau_1,\ldots,\tau_n,\zeta_1,\ldots,\zeta_n) 1_{\tau_n\leq t},
 \enqs
 for $a$ $=$ $(a_0,\ldots,a_n)$ $\in$ $A$ $=$ $A_0\times \ldots \times A_n$,  where $f^0$ is $\Oc(\F)\otimes\Bc(\R^d)\otimes\Bc(A_0)$-measurable, and 
 $f^k$ is $\Oc(\F)\otimes\Bc(\R^d)\otimes\Bc(A_k)\otimes\Bc(\Delta_k)\otimes\Bc(E^k)$-measurable, for  $k$ $=$ $1,\ldots,n$.

 The value function for the stochastic control problem is then defined by:
\beq \label{defvalue}
V_0(x) &=& \sup_{\alpha \in \Ac_{\G}} \E\Big[ \int_0^T f_t(X_t^{x,\alpha},\alpha_t) dt + G_T(X_T^{x,\alpha}) \Big], \;\;\; x \in \R^d.  
\enq

\begin{Remark}
{\rm In the formulation \reff{defvalue} of  our stochastic control problem, there is a change of regimes in the running and terminal gain 
after each default time.  This is in the spirit of the recent concept of forward or progressive utility functions introduced in  \cite{muszar}.  
}
\end{Remark}

\section{$\F$-decomposition of the stochastic control problem}

\setcounter{equation}{0} \setcounter{Assumption}{0}
\setcounter{Theorem}{0} \setcounter{Proposition}{0}
\setcounter{Corollary}{0} \setcounter{Lemma}{0}
\setcounter{Definition}{0} \setcounter{Remark}{0}

In this section, we provide a decomposition of the value function for the stochastic control problem in the $\G$-filtration, defined in \reff{defvalue}, 
that we  formulate in a  backward induction for value functions of stochastic control in the $\F$-filtration.  To alleviate 
 notations, we shall often omit in \reff{Xkcon} the dependence of $X^{k,x}$ on $\alpha^k$ and 
 $(\theta_1,\ldots,\theta_k,e_1,\ldots,e_k)$ when there is no ambiguity.

\begin{Theorem} \label{thm1}
The value function $V_0$   is obtained from  the  backward induction formula: 
\beq
V_n(x,\theta,e) &=& \esssup_{\alpha^n\in\Ac_{\F}^n} 
\E\Big[ \int_{\theta_n}^T f_t^n(X_t^{n,x},\alpha_t^n,\theta,e) \gamma_t(\theta,e) dt \nonumber \\
& & \hspace{2cm}  + \; 
   G_T^n(X_T^{n,x},\theta,e) \gamma_T(\theta,e) \Big| \Fc_{\theta_n} \Big]  \label{relVn}  \\
V_k(x,\theta^{(k)},e^{(k)}) &=&   \esssup_{\alpha^k\in\Ac_{\F}^k}  
\E\Big[ \int_{\theta_k}^T f_t^k(X_t^{k,x},\alpha_t^k,\theta^{(k)},e^{(k)}) \gamma_t^k(\theta^{(k)},e^{(k)}) dt  \nonumber \\
& & \;\;\;\;\; + \; 
 G_T^k(X_T^{k,x},\theta^{(k)},e^{(k)}) \gamma_T^k(\theta^{(k)},e^{(k)})  \nonumber \\
& & \; \;\; + \; \int_{\theta_k}^T \int_E 
V_{k+1}\big(\Gamma_{\theta_{k+1}}^{k+1}(X_{\theta_{k+1}}^{k,x},\alpha^k_{\theta_{k+1}},e_{k+1}),\theta^{(k)},\theta_{k+1},e^{(k)},e_{k+1} \big)\nonumber  \\
& & \hspace{2cm} \eta_1(de_{k+1}) d\theta_{k+1} \Big| \Fc_{\theta_k} \Big], \;\;\; k=0,\ldots,n-1, \label{relVk}
\enq
for all $\theta$ $=$ $(\theta_1,\ldots,\theta_n)$ $\in$ $\Delta_n\cap [0,T]^n$, $e$ $=$ $(e_1,\ldots,e_n)$ $\in$ $E^n$, $x$ $\in$ $\R^d$. 
\end{Theorem}

\begin{Remark} \label{remthm1}
{\rm  Each step in the backward  induction for the determination of the original value function $V_0$ leads to the 
formulation of a family of value functions associated to  standard stochastic control problem in the $\F$-filtration.  
Indeed, at step $n$,  $V_n(x,.)$ is a family of value functions  parametrized by $(\theta_1,\ldots,\theta_n)$ $\in$ $\Delta_n$, 
$(e_1,\ldots,e_n)$ $\in$ $E^n$, and corresponding to the  stochastic control problem  after the last default at time $\theta_n$, 
with a running gain function $f_t^n$ and terminal gain function $G_T^n$ on the controlled state process 
$X^n$ in the $\F$-filtration, and weighted by  the $\Oc(\F)$-measurable process  $\gamma$.   Now, suppose that  at step $k+1$, 
we have determined the family of value functions $V_{k+1}(x,.)$, $(\theta_1,\ldots,\theta_{k+1})$ $\in$ $\Delta_{k+1}$, $(e_1,\ldots,e_{k+1})$ $\in$ 
$E^{k+1}$, and denote by $\hat V_{k+1}$ the map on $\Omega\times\R^d\times A_k\times\Delta_{k+1}\times E^k$:
\beqs
& & \hat V_{k+1}\big(x,a_k,\theta^{(k)}, \theta_{k+1},e^{(k)}\big) \\
&=&   \int_E V_{k+1} \big(\Gamma_{\theta_{k+1}}^{k+1}(x,a_k,e_{k+1}),\theta^{(k)},\theta_{k+1},e^{(k)},e_{k+1}\big) \eta_1(de_{k+1}). 
\enqs
Then,  the family of value functions at step $k$,  representing the value for the stochastic control problem after $k$ defaults,  
is computed from the stochastic control problem in the $\F$-filtration with the running gain function $f_t^k$ and terminal 
gain function $G_T^k$ weighted by the $\Oc(\F)$-measurable random variable $\gamma^k$, and with the running gain function 
$\hat V_{k+1}$: 
\beq
V_k(x) &=&   \esssup_{\alpha^k\in\Ac_{\F}^k}  
\E\Big[ \int_{\theta_k}^T f_t^k(X_t^{k,x},\alpha_t^k) \gamma_t^k dt  + G_T^k(X_T^{k,x}) \gamma_T^k    \nonumber \\
& & \;\;\; \;\;\; + \;   \int_{\theta_k}^T \hat V_{k+1}(X_{\theta_{k+1}}^{k,x}, \alpha^k_{\theta_{k+1}},\theta_{k+1})  \label{Vkdec}
d \theta_{k+1} \Big| \Fc_{\theta_k} \Big].
\enq
Here, we omitted the dependence in $\theta^{(k)}$ $= $$(\theta_1,\ldots,\theta_k)$, $e^{(k)}$ $=$ $(e_1,\ldots,e_k)$ to alleviate notations. 
The two first terms in the rhs of \reff{Vkdec} represent the gain functional when there is no more default after the $k$-th one, while the last term 
represents  the gain in the case when a $k+1$-th default would occur between the last one at time $\tau_k$ $=$ $\theta_k$ and the finite horizon $T$.   Finally,  the decomposition in Theorem \ref{thm1} also shows that an optimal control for the global problem in the $\G$-filtration is obtained by 
a concatenation of optimal controls for each subproblems $V_k$ in the $\F$-filtration. 
}
\end{Remark}

\noindent {\bf Proof of Theorem \ref{thm1}}. 
 
\noindent   Fix $x$ $\in$ $\R^d$,  $\alpha$ $=$ $(\alpha^0,\ldots,\alpha^n)$ $\in$ $\Ac_\G$, and consider the controlled state process 
$X^{x,\alpha}$.  
By definition of $X^{x,\alpha}$ in \reff{dynXcon},   $G_T(.)$ and $f_t(.)$, 
observe that the $\Gc_T$-measurable random variable $G_T(X_T^{x,\alpha})$ is decomposed 
according to the $n+1$-tuple $(G_T^0(\bar X_T^0),\ldots,G_T^n(\bar X_T^n))$, and the $\G$-optional process $f_t(X_t^{x,\alpha},\alpha_t)$ is decomposed  as  $(f_t^0(\bar X_t^0,\alpha_t^0),\ldots,f_t^n(\bar X_t^n,\alpha_t^n))$. 
Let us now define by backward induction the  maps $J_k$, $k$ $=$ $0,\ldots,n$ by
\beq
J_n(x,\theta,e,\alpha) &=& \E \Big[ \int_{\theta_n}^T f_t^n(X_t^{n,x},\alpha_t^n,\theta,e) \gamma_t(\theta,e) dt + G_T^n(X_T^{n,x},\theta,e) 
\gamma_T(\theta,e) \Big| \Fc_{\theta_n} \Big] \nonumber \\
J_k(x,\theta^{(k)},e^{(k)},\alpha) &=&  \E\Big[ \int_{\theta_k}^T f_t^k(X_t^{k,x},\alpha_t^k,\theta^{(k)},e^{(k)}) \gamma_t^k(\theta^{(k)},e^{(k)}) dt  
\nonumber \\
& & \;\;\;\;\;\;\; + \;  G_T^k(X_T^{k,x},\theta^{(k)},e^{(k)}) \gamma_T^k(\theta^{(k)},e^{(k)})  \nonumber \\
& &  \; + \;  \int_{\theta_k}^T \int_E  J_{k+1}\big(\Gamma_{\theta_{k+1}}^{k+1}(X_{\theta_{k+1}}^{k,x},\alpha^k_{\theta_{k+1}},e_{k+1}),\theta^{(k)},
\theta_{k+1},e^{(k)},e_{k+1},\alpha \big) \nonumber \\
& & \hspace{2cm} \eta_1(de_{k+1}) d\theta_{k+1} \Big| \Fc_{\theta_k} \Big], \label{inducVk}
\enq
for any $x$ $\in$ $\R^d$,  $\theta$ $=$ $(\theta_1,\ldots,\theta_n)$ $\in$ $\Delta_n\cap [0,T]^n$, $e$ $=$ $(e_1,\ldots,e_n)$ $\in$ $E^n$,  and 
$\alpha$ $=$ $(\alpha^0,\ldots,\alpha^n)$ $\in$ $\Ac_\F^0\times\ldots\times\Ac_\F^n$.  Let us denote by $\bar J_k(\theta^{(k)},e^{(k)})$ $=$ 
$J_k(\bar X_{\theta_k}^k,\theta^{(k)},e^{(k)},\alpha)$, $k$ $=$ $0,\ldots,n$, and observe by definition of $X^{x,\alpha}$ and $\bar X^k$ in \reff{dynXcon}  that $\bar J_k$ satisfy the backward induction formula: 
\beqs
\bar J_n(\theta,e) &=&  \E \Big[ \int_{\theta_n}^T f_t^n(\bar X_t^{n},\alpha_t^n,\theta,e) \gamma_t(\theta,e) dt + G_T^n(\bar X_T^{n},\theta,e) \gamma_T(\theta,e) \Big| \Fc_{\theta_n} \Big] \\
\bar J_k(\theta^{(k)},e^{(k)}) &=&   \E\Big[ \int_{\theta_k}^T f_t^k(\bar X_t^{k},\alpha_t^k,\theta^{(k)},e^{(k)}) \gamma_t^k(\theta^{(k)},e^{(k)}) dt \\
& & \;\;\;\;\;\;\;  + \;   G_T^k(\bar X_T^{k},\theta^{(k)},e^{(k)}) \gamma_T^k(\theta^{(k)},e^{(k)})  \\
& &   \; + \;  \int_{\theta_k}^T \int_E \bar J_{k+1}(\theta^{(k)},\theta_{k+1},e^{(k)},e_{k+1}) \eta_1(de_{k+1}) d\theta_{k+1} \Big| \Fc_{\theta_k} \Big].
\enqs
Therefore, from Proposition \ref{corollexp}, we have the equality: 
\beq \label{relinter}
\E\Big[ \int_0^T f(X_t^{x,\alpha},\alpha_t) dt + G_T(X_T^{x,\alpha})\Big] \; = \;  \bar J_0 &=& J_0(x,\alpha).  
\enq

Let us now define the value function processes: 
\beq \label{defVk}
V_k(x,\theta^{(k)},e^{(k)}) & := & \esssup_{\alpha\in\Ac_\G} J_k(x,\theta^{(k)},e^{(k)},\alpha),
\enq
for $k$ $=$ $0,\ldots,n$, $x$ $\in$ $\R^d$, and $\theta$ $=$ $(\theta_1,\ldots,\theta_n)$ $\in$ $\Delta_n\cap [0,T]^n$, 
$e$ $=$ $(e_1,\ldots,e_n)$ $\in$ $E^n$.  First, observe that this definition for $k$ $=$ $0$ is consistent with the definition of the value function $V_0$ 
of the stochastic control problem \reff{defvalue} from the relation \reff{relinter}.  Then, it  remains to  prove that the value functions $V_k$ defined 
in \reff{defVk} satisfy the backward induction formula in the assertion of the theorem.  For $k$ $=$ $n$, and since $J_n(x,\theta,e,\alpha)$ depends on $\alpha$ only through its last component $\alpha^n$,  the relation \reff{relVn} holds true.  Next, from the backward induction \reff{inducVk} for 
$J_k$, and the definition of $V_{k+1}$, we have for all $\alpha$ $=$ $(\alpha^0,\ldots,\alpha^n)$ $\in$ $\Ac_\G$:
\beq
J_k(x,\theta^{(k)},e^{(k)},\alpha) & \leq &  \E\Big[ \int_{\theta_k}^T f_t^k(X_t^{k,x},\alpha_t^k,\theta^{(k)},e^{(k)}) \gamma_t^k(\theta^{(k)},e^{(k)}) dt  
\nonumber \\
& & \;\;\;\;\;\;\; + \;  G_T^k(X_T^{k,x},\theta^{(k)},e^{(k)}) \gamma_T^k(\theta^{(k)},e^{(k)})  \nonumber \\
& &  \; + \;  \int_{\theta_k}^T \int_E  V_{k+1}\big(\Gamma_{\theta_{k+1}}^{k+1}(X_{\theta_{k+1}}^{k,x},\alpha^k_{\theta_{k+1}},e_{k+1}),\theta^{(k)},
\theta_{k+1},e^{(k)},e_{k+1} \big)  \nonumber \\
& & \hspace{2cm} \eta_1(de_{k+1}) d\theta_{k+1} \Big| \Fc_{\theta_k} \Big]  \nonumber \\
& \leq & \bar V_k(x,\theta^{(k)},e^{(k)}), \label{JkbarVk}
\enq
where $\bar V_k$ is defined by  the rhs of \reff{relVk}.  By taking the supremum over $\alpha$ in the inequality \reff{JkbarVk}, this shows that 
$V_k$ $\leq$ $\bar V_k$.  Conversely, fix $x$ $\in$ $\R^d$, $\theta$ $=$ $(\theta_1,\ldots,\theta_n)$ $\in$ $\Delta_n\cap [0,T]^n$, 
$e$ $=$ $(e_1,\ldots,e_n)$ $\in$ $E^n$, and let us prove that $V_k(x,\theta^{(k)},e^{(k)})$ $\geq$ $\bar V_k(x,\theta^{(k)},e^{(k)})$.  
Fix an arbitrary $\alpha^k$ $\in$ $\Ac_\F^k$, and the associated controlled process $X^{k,x}$. 
By definition of $V_{k+1}$, for any $\omega$ $\in$ $\Omega$, $\eps$ $>$ $0$, there exists $\alpha^{\omega,\eps}$ $\in$ $\Ac_\G$, which is an 
$\eps$-optimal control for $V_{k+1}(.,\theta^{(k)},e^{(k)})$ at 
$(\omega,\Gamma_{\theta_{k+1}}^{k+1}(X_{\theta_{k+1}}^{k,x},\alpha^k_{\theta_{k+1}},e_{k+1}))$. Recalling that the set of admissible controls is a separable metric space, one can use a measurable selection result (see e.g. \cite{wag80}) to find $\alpha^\eps$ $\in$ $\Ac_\G$ s.t. 
$\alpha^\eps_t(\omega)$ $=$ $\alpha_t^{\omega,\eps}(\omega)$,  $dt\otimes d\P$ a.e., and so
\beqs
& & V_{k+1}\big(\Gamma_{\theta_{k+1}}^{k+1}(X_{\theta_{k+1}}^{k,x},\alpha^k_{\theta_{k+1}},e_{k+1}),\theta^{(k)},\theta_{k+1},e^{(k)},e_{k+1} \big) 
-\eps \\
&\leq &  J_{k+1}\big(\Gamma_{\theta_{k+1}}^{k+1}(X_{\theta_{k+1}}^{k,x},\alpha^k_{\theta_{k+1}},e_{k+1}),\theta^{(k)},
\theta_{k+1},e^{(k)},e_{k+1},\alpha^\eps \big), \;\;\;\;\; a.s. 
\enqs
 Denote by $(\alpha^{\eps,0},\ldots,\alpha^{\eps,n})$ the $n+1$-tuple associated to $\alpha^\eps$ $\in$ $\Ac_\G$, and let us consider the admissible control $\tilde\alpha^\eps$ $=$ $(\alpha^{\eps,0},\ldots,\alpha^k,\alpha^{\eps,k+1},\ldots,\alpha^{\eps,n})$ $\in$ $\Ac_\G$ consisting in  substituting the $k$-th component of $\alpha^\eps$ by $\alpha^k$ $\in$ $\Ac_\F^k$. Since $J_{k+1}(x,\theta,e,\alpha)$ depends on $\alpha$ only through its 
last components $(\alpha^{k+1},\ldots,\alpha^n)$, we have from \reff{inducVk}
\beqs
V_k(x,\theta^{(k)},e^{(k)}) & \geq & J_k(x,\theta^{(k)},e^{(k)},\tilde \alpha^\eps) \\
&=& \E\Big[ \int_{\theta_k}^T f_t^k(X_t^{k,x},\alpha_t^k,\theta^{(k)},e^{(k)}) \gamma_t^k(\theta^{(k)},e^{(k)}) dt  
\nonumber \\
& & \;\;\;\;\;\;\; + \;  G_T^k(X_T^{k,x},\theta^{(k)},e^{(k)}) \gamma_T^k(\theta^{(k)},e^{(k)})  \nonumber \\
& &  \; + \;  \int_{\theta_k}^T \int_E  J_{k+1}\big(\Gamma_{\theta_{k+1}}^{k+1}(X_{\theta_{k+1}}^{k,x},\alpha^k_{\theta_{k+1}},e_{k+1}),\theta^{(k)},
\theta_{k+1},e^{(k)},e_{k+1},\alpha^\eps \big) \nonumber \\
& & \hspace{2cm} \eta_1(de_{k+1}) d\theta_{k+1} \Big| \Fc_{\theta_k} \Big] \\
& \geq & \E\Big[ \int_{\theta_k}^T f_t^k(X_t^{k,x},\alpha_t^k,\theta^{(k)},e^{(k)}) \gamma_t^k(\theta^{(k)},e^{(k)}) dt  
\nonumber \\
& & \;\;\;\;\;\;\; + \;  G_T^k(X_T^{k,x},\theta^{(k)},e^{(k)}) \gamma_T^k(\theta^{(k)},e^{(k)})  \nonumber \\
& &  \; + \;  \int_{\theta_k}^T \int_E  V_{k+1}\big(\Gamma_{\theta_{k+1}}^{k+1}(X_{\theta_{k+1}}^{k,x},\alpha^k_{\theta_{k+1}},e_{k+1}),\theta^{(k)},
\theta_{k+1},e^{(k)},e_{k+1} \big)  \nonumber \\
& & \hspace{2cm} \eta_1(de_{k+1}) d\theta_{k+1} \Big| \Fc_{\theta_k} \Big]  - \eps.
\enqs
From the arbitrariness of $\alpha^k$ $\in$ $\Ac_\F^k$ and $\eps$ $>$ $0$, we obtain the required inequality: $V_k(x,\theta^{(k)},e^{(k)})$ $\geq$ 
$\bar V_k(x,\theta^{(k)},e^{(k)})$, and the proof is complete. 
\ep

\section{The case of enlarged filtration with multiple random times } \label{secmultiple}

\setcounter{equation}{0} \setcounter{Assumption}{0}
\setcounter{Theorem}{0} \setcounter{Proposition}{0}
\setcounter{Corollary}{0} \setcounter{Lemma}{0}
\setcounter{Definition}{0} \setcounter{Remark}{0}

In this section, we consider the case where the random times are not  assumed to be  ordered.  In other words, this means that one has 
access to the default times themselves with their indexes, and not only to the ranked default times.   
This general case can actually be derived from the case of successive random times associated with suitable auxiliary marks.  Let us 
consider the progressive enlargement of filtration  from $\F$ to $\G$ with multiple random times $(\tau_1,\ldots,\tau_n)$ associated with the 
marks $(\zeta_1,\ldots,\zeta_n)$. Denote by $\hat\tau_1\leq\ldots\leq\hat\tau_n$ the corresponding ranked  times, and by 
$\iota_i$ the index mark (valued in $\{1,\ldots,n\}$) of  the $i$-th order statistic of $(\tau_1,\ldots,\tau_n)$ for $i$ $=$ $1,\ldots,n$.  Then, it is clear that 
the progressive enlargement of filtration of $\F$ with the successive random times $(\hat\tau_1,\ldots,\hat\tau_n)$ together with the marks 
$(\iota_1,\zeta_{\iota_1},\ldots,\iota_n,\zeta_{\iota_n})$ leads to the filtration $\G$, so that one can apply the results of the previous sections.  
For simplicity of notations, we shall  focus on the case of two random times $\tau_1$ and $\tau_2$,   
associated to the marks $\zeta_1$ and $\zeta_2$ valued in $E$ Borel space of $\R^m$.  

The decomposition of optional and predictable process with respect to this progressive enlargement of filtration is given 
by the following lemma, which is derived from Lemma \ref{lemdecG}, with the specific feature that we have also to take into account 
the index of the order statistic in $(\tau_1,\tau_2)$.

\begin{Lemma}
Any $\G$-optional (resp. predictable) process $Y$ $=$ $(Y_t)_{t\geq 0}$ is represented as
\beqs
Y_t &=& Y_t^0 1_{t<\hat\tau_1} + Y_t^{1,1}(\tau_1,\zeta_1) 1_{\tau_1\leq t<\tau_2} +  Y_t^{1,2}(\tau_2,\zeta_2)  1_{\tau_2\leq t<\tau_1} 
+ Y_t^2(\tau_1,\tau_2,\zeta_1,\zeta_2) 1_{t \geq \hat\tau_2},   \\
& (\mbox{resp.} \; =  & Y_t^0 1_{t\leq\hat\tau_1} + Y_t^{1,1}(\tau_1,\zeta_1) 1_{\tau_1< t\leq \tau_2} 
+  Y_t^{1,2}(\tau_2,\zeta_2)  1_{\tau_2< t\leq \tau_1} 
+ Y_t^2(\tau_1,\tau_2,\zeta_1,\zeta_2) 1_{t > \hat\tau_2}), 
\enqs
for all $t$ $\geq$ $0$, where $Y^0$ $\in$ $\Oc_{\F}$ (resp. $\Pc_\F$), 
$Y^{1,1}$, $Y^{1,2}$ $\in$ $\Oc^1_{\F}(\R_+,E)$ (resp. $\Pc_\F^1(\R_+,E)$), and $Y^2$ $\in$ $\Oc^2_{\F}(\R_+^2,E^2)$ (resp. $\Pc_\F^2(\R_+^2,E^2)$).  
\end{Lemma}

 Any $Y$ $\in$ $\Oc_{\G}$ (resp. $\Pc_\G$) can then be identified with a quadruple $(Y^0,Y^{1,1},Y^{1,2},Y^2)$ $\in$ 
 $\Oc_{\F}\times\Oc^1_{\F}(\R_+,E)\times\Oc^1_{\F}(\R_+,E)\times\Oc^2_{\F}(\R_+^2,E^2)$ (resp. 
 $\Pc_{\F}\times\Pc^1_{\F}(\R_+,E)\times\Pc^1_{\F}(\R_+,E)\times\Pc^2_{\F}(\R_+^2,E^2)$).

 \vspace{2mm}

 Similarly as in Section 1, we now make a density hypothesis on the conditional distribution of $(\tau_1,\tau_2,\zeta_1,\zeta_2)$ given the 
 reference information. We assume that  there exists a $\Oc(\F)\otimes\Bc(\R_+^2)\otimes\Bc(E^2)$-measurable 
 map $(t,\omega,\theta_1,\theta_2,e_1,e_2)$ $\rightarrow$ $\gamma_t(\omega,\theta_1,\theta_2,e_1,e_2)$ such that
 \beqs
{\bf (DH)} \hspace{7mm} 
\P\big[ (\tau_1,\tau_2,\zeta_1,\zeta_2) \in d\theta  de | \Fc_t \big] 
&=& \gamma_t(\theta_1,\theta_2,e_1,e_2)  d\theta_1d\theta_2 \eta(de_1)\eta(de_2), \;\;\; a.s. 
\enqs
where $\eta$ is a nonnegative measure on $\Bc(E)$.

We next introduce some useful notations. We denote by $\gamma_0$ the $\F$-optional process defined by
\beqs
\gamma_t^0 &=& \P[\tau_1 > t |\Fc_t] \; = \; \int_{E^2} \int_{[t,\infty)^2} \gamma_t(\theta_1,\theta_2,e_1,e_2) d\theta_1d\theta_2
 \eta(de_1)\eta(de_2),
\enqs
and we denote by $(t,\omega,\theta_1,e_1)$ $\rightarrow$ $\gamma_t^{1,1}(\theta_1,e_1)$,  
and $(t,\omega,\theta_2,e_2)$ $\rightarrow$ $\gamma_t^{1,2}(\theta_2,e_2)$, $t$ $\geq$ $0$, the 
$\Oc(\F)\otimes\Bc(\R_+)\otimes\Bc(E)$-measurable maps defined by
\beqs
\gamma_t^{1,1}(\theta_1,e_1) &= &  \int_E\int_t^\infty \gamma_t(\theta_1,\theta_2,e_1,e_2) d\theta_2 \eta(de_2), \\
\gamma_t^{1,2}(\theta_2,e_2) &=& \int_E\int_t^\infty \gamma_t(\theta_1,\theta_2,e_1,e_2) d\theta_1 \eta(de_1),
\enqs
so that
\beqs
\P[ \tau_2>t | \Fc_t] \; = \; \int_E \int_0^\infty \gamma_t^{1,1}(\theta_1,e_1) d\theta_1\eta(de_1), & & 
\P[ \tau_1>t | \Fc_t] \; = \; \int_E \int_0^\infty \gamma_t^{1,2}(\theta_2,e_2) d\theta_2\eta(de_2).
\enqs

 \vspace{2mm}
 
 The next result, which is analog to Proposition \ref{corollexp},  provides a backward induction formula involving $\F$-expectations for the computation of expectation functionals of  $\G$-optional processes.

 \begin{Proposition}
 Let $Y$ $=$ $(Y^0,Y^{1,1},Y^{1,2},Y^2)$ and $Z$ $=$ $(Z^0,Z^{1,1},Z^{1,2},Z^2)$ be two nonnegative (or bounded) $\G$-optional processes, 
 and fix $T$ $\in$ $(0,\infty)$. 
 
 \noindent The expectation $\E[\int_0^T Y_t dt + Z_T]$ can be computed in a backward induction as
 \beqs
\E\Big[\int_0^T Y_t dt + Z_T\Big] &=&  J_0
\enqs
where the  $(J_0,J_{1,1},J_{1,2},J_2)$ are given by
 \beqs
 J_{2}(\theta_1,\theta_2,e_1,e_2) &=& \E \Big[ \int_{\theta_1\vee\theta_2}^T Y_t^2  \gamma_t(\theta_1,\theta_2,e_1,e_2) dt    +    
 Z_T^2 \gamma_T(\theta_1,\theta_2,e_1,e_2) \Big| \Fc_{\theta_1\vee\theta_2} \Big]  \\
 J_{1,1}(\theta_1,e_1) &=& \E \Big[ \int_{\theta_1}^T Y_t^{1,1} \gamma_t^{1,1} (\theta_1,e_1)  dt + Z_T^{1,1} \gamma_T^{1,1} (\theta_1,e_1)  \\
 & & \;\;\;  + \;  \int_E  \int_{\theta_1}^T J_{2}(\theta_1,\theta_2,e_1,e_2) d \theta_2\eta(de_2)  \Big| \Fc_{\theta_1} \Big] \\
 J_{1,2}(\theta_2,e_2) &=& \E \Big[ \int_{\theta_2}^T Y_t^{1,2} \gamma_t^{1,2} (\theta_2,e_2)  dt + Z_T^{1,2} \gamma_T^{1,2} (\theta_2,e_2)  \\
 & & \;\;\;  + \;  \int_E  \int_{\theta_1}^T J_{2}(\theta_1,\theta_2,e_1,e_2) d \theta_2\eta(de_2)  \Big| \Fc_{\theta_2} \Big] \\
 J_{0} &=& \E\Big[ \int_0^T Y_t^0 \gamma_t^0 dt + Z_T^0 \gamma_T^0  \\
 & & \;\;\; + \;   \int_E \int_0^T J_{1,1}(\theta_1,e_1) d\theta_1\eta(de_1) +  \int_E \int_0^T J_{1,2}(\theta_2,e_2) d\theta_2\eta(de_2) \Big].  
 \enqs
 \end{Proposition}

\vspace{2mm}

Let us now formulate the general stochastic control problem in this framework.   

\noindent A control is a  $\G$-predictable process $\alpha$ $=$ $(\alpha^0,\alpha^{1,1},\alpha^{1,2},\alpha^2)$ $\in$ 
$\Pc_{\F}\times\Pc^1_{\F}(\R_+,E)\times\Pc^1_{\F}(\R_+,E)\times\Pc^2_{\F}(\R_+^2,E^2)$, where $\alpha^0$, $\alpha^{1,1}$, $\alpha^{1,2}$ and  
$\alpha^2$  are valued respectively in $A_0$,  $A_{1,1}$, $A_{1,2}$ and $A_2$, Borel sets of some Euclidian space. 
We denote by $A$ $=$ $A_0\times A_{1,1}\times A_{1,2}\times A_2$, and by 
$\Ac_{\G}$ the set of admissible control processes, which is a product space 
$\Ac_{\F}^0\times\Ac_{\F}^{1,1}\times\Ac_{\F}^{1,2}\times\Ac_{\F}^2$, where  $\Ac_{\F}^0$, $\Ac_{\F}^{1,1}$, $\Ac_{\F}^{1,2}$ and $\Ac_{\F}^2$ 
are some separable metric spaces respectively in  $\Pc_{\F}(A_0)$, $\Pc^1_{\F}(\R_+,E;A_{1,1})$, 
$\Pc^1_{\F}(\R_+,E;A_{1,2})$ and  $\Pc^2_{\F}(\R_+^2,E^2;A_2)$. 
 
We are next given a collection of measurable mappings:
\beqs
(x,\alpha^0) \in \R^d\times\Ac_{\F}^0 & \longmapsto & X^{0,x,\alpha^0} \; \in \;  \Oc_{\F} \\
(x,\alpha^{1,1}) \in \R^d\times\Ac_{\F}^{1,1} & \longmapsto & X^{1,1,x,\alpha^{1,1}} \; \in \; \Oc^{1}_{\F}(\R_+,E) \\
(x,\alpha^{1,2}) \in \R^d\times\Ac_{\F}^{1,2} & \longmapsto & X^{1,2,x,\alpha^{1,2}} \; \in \; \Oc^{1}_{\F}(\R_+,E) \\
(x,\alpha^2) \in \R^d\times\Ac_{\F}^2 & \longmapsto & X^{2,x,\alpha^2} \; \in \; \Oc^2_{\F}(\R_+^2,E^2),
\enqs
such that 
we have the initial data
\beqs
X_0^{0,x,\alpha^0} &=& x,  \;\;\;  \forall x \in \R^d, \\
X_{\theta_1}^{1,1,\xi,\alpha^{1,1}}(\theta_1,e_1) &=& \xi, \;\;\; \forall \xi \; \Fc_{\theta_1}-\mbox{measurable}, \\
X_{\theta_2}^{1,2,\xi,\alpha^{1,2}}(\theta_2,e_2) &=& \xi, \;\;\; \forall \xi \; \Fc_{\theta_2}-\mbox{measurable},  \\
X_{\theta_1\vee\theta_2}^{2,\xi,\alpha^2}(\theta_1,\theta_2,e_1,e_2) &=&  \xi, \;\;\;  \forall \xi \; \Fc_{\theta_1\vee\theta_2}-\mbox{measurable}.
\enqs

We are also given a collection of maps $\Gamma^{1,1}$, $\Gamma^{1,2}$,  on  $\R_+\times\Omega\times\R^d\times A_0\times E$, 
$\Gamma^{2,1}$  on $\R_+\times\Omega\times\R^d\times A_{1,1}\times E$ and 
$\Gamma^{2,2}$ on  $\R_+\times\Omega\times\R^d\times A_{1,2}\times E$ such that  
\beqs
(t,\omega,x,a,e) & \mapsto & \Gamma^{1,1}_t(\omega,x,a,e), \; \Gamma^{1,2}_t(\omega,x,a,e)  \\
& & \;\;\;\; \mbox{ are  } \;\;
\Pc(\F)\otimes\Bc(\R^d)\otimes\Bc(A_{0})\otimes\Bc(E)-\mbox{measurable} \\
(t,\omega,x,a,e) & \mapsto & \Gamma^{2,1}_t(\omega,x,a,e)  \; \mbox{ is } \; 
\Pc(\F)\otimes\Bc(\R^d)\otimes\Bc(A_{1,1})\otimes\Bc(E)-\mbox{measurable} \\
(t,\omega,x,a,e) & \mapsto & \Gamma^{2,2}_t(\omega,x,a,e)  \; \mbox{ is } \; 
\Pc(\F)\otimes\Bc(\R^d)\otimes\Bc(A_{1,2})\otimes\Bc(E)-\mbox{measurable}
\enqs

The controlled state process is then given by the mapping
\beqs
(x,\alpha) \in  \R^d\times\Ac_{\G} & \longmapsto & X^{x,\alpha} \; \in \; \Oc_{\G},
\enqs
where for $\alpha$ $=$ $(\alpha^0,\alpha^{1,1},\alpha^{1,2},\alpha^2)$, $X^{x,\alpha}$  is the  process equal to
\beqs
X_t^{x,\alpha} &=& \bar X_t^0 1_{t<\hat\tau_1} + \bar X_t^{1,1}(\tau_1,\zeta_1) 1_{\tau_1\leq t<\tau_2} +  
\bar X_t^{1,2}(\tau_2,\zeta_2) 1_{\tau_2\leq t<\tau_1} + \bar X_t^2(\tau_1,\tau_2,\zeta_1,\zeta_2) 1_{t\geq\hat\tau_2},
\enqs
with $(\bar X^0,\bar X^{1,1},\bar X^{1,2},\bar X^2)$ $\in$ 
$\Oc_{\F}\times\Oc^1_{\F}(\R_+,E)\times\Oc^1_{\F}(\R_+,E)\times\Oc^2_{\F}(\R_+^2,E^2)$ given by
\beqs
\bar X_0 &=& X^{0,x,\alpha^0} \\
\bar X^{1,1}(\theta_1,e_1) &=& X^{1,1,\Gamma^{1,1}_{\theta_1}(\bar X^0_{\theta_1},\alpha^0_{\theta_1},e_1),\alpha^{1,1}}(\theta_1,e_1) \\
\bar X^{1,2}(\theta_2,e_2) &=& X^{1,2,\Gamma^{1,2}_{\theta_2}(\bar X^0_{\theta_2},\alpha^0_{\theta_2},e_2),\alpha^{1,2}}(\theta_2,e_2) \\
\bar X^2(\theta_1,\theta_2,e_1,e_2) &=& \left\{ 
								\begin{array}{ll}
X^{2,\Gamma^{2,2}_{\theta_2}(\bar X_{\theta_2}^{1,1},\alpha^{1,1}_{\theta_2},e_2),\alpha^2}(\theta_1,\theta_2,e_1,e_2) & \mbox{ if } \theta_1 \leq \theta_2 \\
X^{2,\Gamma^{2,1}_{\theta_1}(\bar X_{\theta_1}^{1,2},\alpha^{1,2}_{\theta_1},e_1),\alpha^2}(\theta_1,\theta_2,e_1,e_2) & \mbox{ if } \theta_2 < \theta_1. 
								\end{array}
								\right. 
\enqs
The interpretation is the following: $X^{0}$ is the controlled state process before any default, $X^{1,1}$ (resp. $X^{1,2}$)  
is the controlled state process  between $\tau_1$ and $\tau_2$ (resp. between $\tau_2$ and $\tau_1$)  
if  the default of index $1$ (resp. index $2$) occurs  first,  and $X^2$ is the controlled state process after both defaults. Moreover,  $\Gamma^{1,1}$ 
(resp. $\Gamma^{1,2}$) represents the  jump of  $X^0$ at $\tau_1$ (resp. $\tau_2$)  if  the default of index  $1$ (resp. index $2$) occurs first, and 
$\Gamma^{2,2}$ (resp. $\Gamma^{2,1}$) represents the jump of $X^{1,1}$ (resp. $X^{1,2}$) at $\tau_2$ (resp. $\tau_1$)  when the default of index $2$ (resp. index $1$) occurs in second after index $1$ (resp. index $2$).

\vspace{2mm}

For a fixed finite horizon $T$ $<$ $\infty$, we are given a nonnegative map $G_T$ on $\Omega\times\R^d$ such that $(\omega,x)$ $\mapsto$ 
$G_T(\omega,x)$ is $\Gc_T\otimes\Bc(\R^d)$-measurable, thus in the form
\beqs
G_T(x) &=& G_T^0(x) 1_{T < \hat\tau_1} +  G_T^{1,1}(x,\tau_1,\zeta_1) 1_{\tau_1\leq T<\tau_2} + G_T^{1,2}(x,\tau_2,\zeta_2) 1_{\tau_2\leq T<\tau_1} 
\\
& & \;\;\; + \;   G_T^2(x,\tau_1,\tau_2,\zeta_1,\zeta_2) 1_{\hat\tau_2\leq T},
\enqs
where $G_T^0$ is $\Fc_T\otimes\Bc(\R^d)$-measurable, $G_T^{1,1}$, $G_T^{1,2}$ are $\Fc_T\otimes\Bc(\R^d)\otimes\Bc(\R_+)\otimes\Bc(E)$-measurable, and $G_T^2$ is $\Fc_T\otimes\Bc(\R^d)\otimes\Bc(\R_+^2)\otimes\Bc(E^2)$-measurable. 
The running gain function is given by a nonnegative map $f$ on $\Omega\times\R^d\times A$ such that 
 $(t,\omega,x,a)$ $\mapsto$ $f_t(\omega,x,a)$ is $\Oc(\G)\otimes\Bc(\R^d)\otimes\Bc(A)$-measurable, and which may be decomposed as 
 \beqs
 f_t(x,a) &=& f_t^0(x,a_0) 1_{t < \hat\tau_1} +  f_t^{1,1}(x,a_{1,1},\tau_1,\zeta_1) 1_{\tau_1\leq t<\tau_2} + 
 f_t^{1,2}(x,a_{1,2},\tau_2,\zeta_2) 1_{\tau_2\leq t<\tau_1} \\
 & & \;\;\; + \;   f_t^2(x,a_2,\tau_1,\tau_2,\zeta_1,\zeta_2) 1_{\hat\tau_2\leq T},
 \enqs
 for $a$ $=$ $(a_0,a_{1,1},a_{1,2},a_2)$ $\in$ $A$ $=$ $A_0\times A_{1,1}\times A_{1,2}\times A_2$,  
 where $f^0$ is $\Oc(\F)\otimes\Bc(\R^d)\otimes\Bc(A_0)$-measurable, and 
 $f^{1,1}$ is  $\Oc(\F)\otimes\Bc(\R^d)\otimes\Bc(A_{1,1})\otimes\Bc(\R_+)\otimes\Bc(E)$-measurable, 
 $f^{1,2}$ is  $\Oc(\F)\otimes\Bc(\R^d)\otimes\Bc(A_{1,2})\otimes\Bc(\R_+)\otimes\Bc(E)$-measurable
 and  $f^{2}$ is $\Oc(\F)\otimes\Bc(\R^d)\otimes\Bc(A_2)\otimes\Bc(\R_+^2)\otimes\Bc(E^2)$-measurable.

The value function for the stochastic control problem is  then defined by
\beqs
V_0(x) &=& \sup_{\alpha\in\Ac_{\G}} \E\Big[ \int_0^T f_t(X_t^{x,\alpha},\alpha_t) dt + G_T(X_T^{x,\alpha}) \Big], \;\;\; x \in \R^d. 
\enqs

The main result of this section provides a decomposition of the value function in the reference filtration, which is analog to the decomposition in Theorem \ref{thm1}.   To alleviate the notations, we omit the 
dependence of the state process in the controls and in the parameters $\theta,e$, when there is no ambiguity.

\begin{Theorem} \label{thm2}
The value function $V_0$ is obtained from the backward induction formula
\beqs
V_2(x,\theta_1,\theta_2,e_1,e_2) &=& \esssup_{\alpha^2\in\Ac^2_{\F}} \E\Big[ 
\int_{\theta_1\vee\theta_2}^T f_t^2(X_t^ {2,x},\alpha_t^2,\theta_1,\theta_2,e_1,e_2) \gamma_t(\theta_1,\theta_2,e_1,e_2) dt \\
& & \hspace{2cm} + \;  G_T^2(X_T^{2,x},\theta_1,\theta_2,e_1,e_2) 
\gamma_T(\theta_1,\theta_2,e_1,e_2) \Big| 
\Fc_{\theta_1\vee\theta_2} \Big] \\
V_{1,1}(x,\theta_1,e_1) &=& \esssup_{\alpha^{1,1}\in\Ac^{1,1}_{\F}} \E \Big[ 
\int_{\theta_1}^T f_t^{1,1}(X_t^ {1,1,x},\alpha_t^{1,1},\theta_1,e_1) \gamma_t^{1,1}(\theta_1,e_1) dt  \\
& &  \hspace{2cm} + \;  G_T^{1,1}(X_T^{1,1,x},\theta_1,e_1) \gamma_T^{1,1}(\theta_1,e_1) \\
& & \;\;\;\; + \; \int_{\theta_1}^T \int_E V_2\big(\Gamma^{2,2}_{\theta_2}(X_{\theta_2}^{1,1,x},\alpha^{1,1}_{\theta_2},e_2),\theta_1,\theta_2,e_1,e_2
\big) \eta(de_2) d \theta_2 \Big| \Fc_{\theta_1} \Big]  \\
V_{1,2}(x,\theta_2,e_2) &=& \esssup_{\alpha^{1,2}\in\Ac^{1,2}_{\F}} \E \Big[ 
\int_{\theta_2}^T f_t^{1,2}(X_t^ {1,2,x},\alpha_t^{1,2},\theta_2,e_2) \gamma_t^{1,2}(\theta_2,e_2) dt  \\
& &  \hspace{2cm} + \; G_T^{1,2}(X_T^{1,2,x},\theta_2,e_2) \gamma_T^{1,2}(\theta_2,e_2) \\
& & \;\;\;\; + \; \int_{\theta_1}^T \int_E V_2\big(\Gamma^{2,1}_{\theta_1}(X_{\theta_1}^{1,2,x},\alpha^{1,2}_{\theta_1},e_1),\theta_1,\theta_2,e_1,e_2\big) \eta(de_1) d \theta_1 \Big| \Fc_{\theta_2} \Big] \\
V_0(x) &=& \sup_{\alpha^0\in\Ac_{\F}^0} \E\Big[ \int_0^T f_t^0(X_t ^{0,x},\alpha_t^0) \gamma_t^0  dt  +  G_T^0(X_T^{0,x}) \gamma_T^0  \\
& &\hspace{2cm} + \;  \int_0^T \int_E V_{1,1}\big(\Gamma^{1,1}_{\theta_1}(X^{0,x}_{\theta_1},\alpha^0_{\theta_1},e_1),\theta_1,e_1\big) \eta(de_1) d\theta_1 \\
& & \hspace{2cm} + \;  \int_0^T \int_E V_{1,2}\big(\Gamma^{1,2}_{\theta_2}(X^{0,x}_{\theta_2},\alpha^0_{\theta_2},e_2),\theta_2,e_2\big) \eta(de_2) d\theta_2\Big],
\enqs
for all $(\theta_1,\theta_2)$ $\in$ $[0,T]^2$, $(e_1,e_2)$ $\in$ $E^2$. 
\end{Theorem}

\begin{Remark}
{\rm  As mentioned in Remark \ref{remthm1}, the value functions $V_2$, $V_{1,1}$ and $V_{1,2}$ correspond to standard stochastic control problem in the $\F$-filtration. This is also the case for $V_0$ in the decomposition formula of Theorem \ref{thm2}. Indeed,  denote by  $V_1$ the map 
on $\Omega\times[0,T]\times\R^d\times A_0$:  
\beqs
V_1(x,\theta,a_0) &=& \int_E V_{1,1}(\Gamma_\theta^{1,1}(x,a_0,e),\theta,e) + V_{1,2}(\Gamma_\theta^{1,2}(x,a_0,e),\theta,e) \; \eta(de). 
\enqs
Then, $V_0$ is computed from the stochastic control problem in the $\F$-filtration with the terminal gain function $G_T^0$ weighted by the $\Fc_T$-measurable random variable $\gamma_T^0$, and with the running gain functions $f^0\gamma^0$ and $V_1$: 
\beqs
V_0(x) &=& \sup_{\alpha^0\in\Ac_{\F}} \E\Big[ G_T^0(X_T^{0,x}) \gamma_T^0 + \int_0^T f_t^0(X_t ^{0,x},\alpha_t^0) \gamma_t^0 
+ V_1(X_t^{0,x},t,\alpha_t^0) dt \Big]. 
\enqs
}
\end{Remark}

\section{Applications in mathematical finance}

\setcounter{equation}{0} \setcounter{Assumption}{0}
\setcounter{Theorem}{0} \setcounter{Proposition}{0}
\setcounter{Corollary}{0} \setcounter{Lemma}{0}
\setcounter{Definition}{0} \setcounter{Remark}{0}

\subsection{Indifference pricing of defaultable claims}

We consider a stock subject to a single counterparty default at a random time $\tau$, which induces a jump of random relative size $\zeta$ valued in 
$E$ $\subset$ $(-1,\infty)$.   The price process of the stock is described by
\beqs
S_t &=& S_t^0 1_{t< \tau} + S_t^1(\tau,\zeta) 1_{t\geq \tau},
\enqs
where $S^0$ is governed by 
\beqs
dS_t^0 &=& S_t^0 \big( b_t^0 dt + \sigma_t^0 dW_t \big),
\enqs
and the indexed process $S^1(\theta,e)$, $(\theta,e)$ $\in$ $\R_+\times E$ is given by
\beqs
dS_t^1(\theta,e) &=& S_t^1(\theta,e) \big(  b_t^1(\theta,e) dt + \sigma_t^1(\theta,e) dW_t \big), \;\; t \geq \theta, \\
 S_\theta^1(\theta,e) & = &  S_\theta^0.(1+e). 
\enqs
Here $W$ is a $(\P,\F)$-Brownian motion, $b^0$, $\sigma^0$ $>$ $0$ are $\F$-adapted processes, $b^1$, $\sigma^1$ $>$ $0$  $\in$ 
$\Oc_{\F}^1(\R_+,E)$.  The market information is represented by  the progressive enlarged filtration $\G$ $=$ $\F$ $\vee$ $\D$, with 
$\D$ $=$ $(\Dc_t)_{t\geq 0}$, $\Dc_t$ $=$ $\cap_{\eps >0}\{\sigma(\zeta 1_{\tau\leq s}, s\leq t+\eps)\vee\sigma(1_{\tau\leq s},s\leq t+\eps)\}$.  
An investor can trade in a riskless bond with zero interest rate, and in 
the defaultable  stock. Her trading strategy is  a $\G$-predictable process $\alpha$ $=$ $(\alpha^0,\alpha^1)$ $\in$ 
$\Pc_{\F}\times\Pc^1_{\F}(\R_+,E)$ representing the amount traded in the stock.  We allow constraints on  trading strategy by considering  
 closed sets  $A_0$  and $A_1$ in which  the controls $\alpha^0$ and $\alpha^1$ take values. Notice also that  $A_0$ and $A_1$  may 
 differ. The controlled wealth process of the investor is then given by
\beq \label{Xwealth}
X_t &=& X_t^0 1_{t<\tau} + X_t^1(\tau,\zeta) 1_{t\geq\tau},
\enq
where $X^0$ is the wealth process before the default, and governed by 
\beqs
dX_t^0 &=& \alpha_t^0 \frac{dS_t^0}{S_t^0} \; = \; \alpha_t^0( b_t^0 dt + \sigma_t^0 dW_t),
\enqs
and $X^1(\theta,e)$ is the wealth indexed process after-default, governed by
\beqs
dX_t^1(\theta,e) &=& \alpha_t^1(\theta,e) \frac{dS_t^1(\theta,e)}{S_t^1(\theta,e)}  \; = \; 
\alpha_t^1(\theta,e) \big(  b_t^1(\theta,e) dt + \sigma_t^1(\theta,e) dW_t \big), \;\;\; t \geq \theta \\
X_\theta^1(\theta,e) &=& X_\theta^0 + \alpha_\theta^0 e. 
\enqs

Let us now consider a defaultable contingent claim with payoff  at maturity $T$ given by
\beqs
H_T &=& H_T^0 1_{T<\tau} + H_T^1(\tau,\zeta) 1_{\tau\leq T},
\enqs
where $H_T^0$ is a bounded $\Fc_T$-measurable random variable, and $H_T^1(,)$ is  a bounded 
$\Fc_T\otimes\Bc(\R_+)\otimes\Bc(E)$-measurable map. We use  the  popular indifference pricing criterion for valuing this defaultable 
claim. We are then given an exponential utility function $U$ on $\R$, i.e. 
\beqs
U(x) &=& - \exp(- p x), \;\; x \in \R, 
\enqs
for some $p$ $>$ $0$, and we consider the optimal investment problem for an agent delivering the defaultable claim at maturity $T$: 
\beq \label{V0H}
V_0^H(x) &=& \sup_{\alpha \in \Ac_{\G}} \E \big[ U(X_T^{x,\alpha} - H_T) \big].
\enq
Here $X^{x,\alpha}$ is the wealth process in \reff{Xwealth} controlled by the trading strategy $\alpha$, and starting from $x$ at time $0$. 
We denote by $V_0$ the value function for the optimal investment problem without the defaultable claim, i.e. when $H_T$ $=$ $0$ 
in \reff{V0H}, and the indifference price for $H_T$ is the amount of initial capital such that the investor is indifferent between holding or not the defaultable claim. It is then defined as the unique number $\pi$ such that 
\beqs
V_0^H(x+\pi) &=& V_0(x). 
\enqs

A similar problem (without unpredictable mark $\zeta$) was recently considered in \cite{limque09} and \cite{ankblaeyr09} by using a 
global $\G$-filtration approach under {\bf (H)} hypothesis, see also \cite{penxu09}. The paper \cite{jiapha09} studied an optimal investment problem with power utility functions under a single counterparty default by using a density approach for decomposing the problem in  the $\F$-filtration.  
We  follow this methodology and solve the stochastic control 
problem \reff{V0H} by applying the $\F$-decomposition method. From Theorem \ref{thm1}, the value function $V_0^H$ 
is obtained in two steps via the resolution of the after-default  problem
\beq \label{V1H}
V_1^H(x,\theta,e) &=& \esssup_{\alpha^1 \in \Ac_{\F}^1} \E \Big[ U\big(X_T^{1,x}(\theta,e) - H_T^1(\theta,e)\big)\gamma_T(\theta,e) \Big| \Fc_\theta \Big],
\enq
and then via the resolution of the before-default problem
\beq \label{V0Hdec} 
V_0^H(x) &=& \sup_{\alpha^0\in\Ac_{\F}^0} \E \Big[ U(X_T^{0,x} - H_T^0) \gamma_T^0 + \int_0^T \int_E V_1^H(X_\theta^{0,x}+\alpha_\theta^0 e,\theta,e)
\eta(de) d \theta \Big]. 
\enq

\vspace{1mm}

\noindent $\bullet$ {\it Solution to the after-default problem.}

\noindent For fixed $(\theta,e)$ $\in$ $[0,T]\times E$, 
problem \reff{V1H} is a classical  utility maximization problem with random endowment in the complete market model after default 
described by the indexed price process $S^1(\theta,e)$.  Indeed, notice that we can remove the positive term $\gamma_T(\theta,e)$ in \reff{V1H} by 
defining the ``modified claim" $\tilde H_T^1(\theta,e)$ $=$ $H_T^1(\theta,e) + \frac{1}{p} \ln \gamma_T(\theta,e)$ so that  
\beq \label{V1Hnew}
V_1^H(x,\theta,e) &=& \esssup_{\alpha^1 \in \Ac_{\F}^1} 
\E \Big[ U\big(X_T^{1,x}(\theta,e) -  \tilde H_T^1(\theta,e)\big) \Big| \Fc_\theta \Big]. 
\enq
This problem was addressed by several methods in the literature, and we know from  dynamic programming  and  BSDE methods 
(see  \cite{elkrou00}  or  \cite{huimkmul05})) that
\beqs
V_1^H(x,\theta,e) &=&  U \big( x- Y_\theta^{1,H}(\theta,e) \big)
\enqs
 where $Y^{1,H}(\theta,e)$ is the unique bounded solution to the BSDE
 \beqs
 Y_t^{1,H}(\theta,e) &=& H_T^1(\theta,e) + \frac{1}{p} \ln \gamma_T(\theta,e)  + \int_t^T f^{1}(r,Z_r^{1,H},\theta,e) dr - \int_t^T Z_r^{1,H} dW_r
 \enqs 
 and the generator $f^{1}$ is the $\Pc(\F)\otimes\Bc(\R)\otimes\Bc(\R_+)\otimes\Bc(E)$-measurable map defined by
 \beqs
 f^{1}(t,z,\theta,e) &=&   -   \frac{b_t^1(\theta,e)}{\sigma_t^1(\theta,e)}  z  - 
 \frac{1}{2p}   \Big(\frac{b_t^1(\theta,e)}{\sigma_t^1(\theta,e)} \Big)^2 +  \frac{p}{2} \inf_{a\in A_1} 
 \Big| \Big( z+ \frac{1}{p} \frac{b_t^1(\theta,e)}{\sigma_t^1(\theta,e)} \Big) -  a \sigma_t^1(\theta,e) \Big|^2. 
 \enqs


\noindent $\bullet$ {\it Global solution}

\noindent The global solution is finally obtained from the resolution of the before-default problem, which is then reduced to 
\beqs
V_0^H(x) &=& \sup_{\alpha^0\in\Ac_{\F}^0} \E \Big[ U(X_T^{0,x} - H_T^0) \gamma_T^0  
+ \int_0^T \int_E  U(X_\theta^{0,x} +  \alpha_\theta^0 e - Y^{1,H}_\theta(e) ) \eta(de) d \theta \Big]. 
\enqs
From the additive dependence of the wealth process $X^{0,x}$ in function of $x$, and the exponential form of the utility function $U$, we know that 
the value function $V_0^H$ is in the form 
\beqs
V_0^H(x) &=& U(x - Y_0^{0,H}),
\enqs
for some quantity $Y_0^{0,H}$ independent of $x$, and which may be characterized by  dynamic programming methods in the $\F$-filtration.  
This can be achieved either via PDE methods in a Markovian setting, or via BSDE methods in the general case.  The BSDE associated to 
$Y^{0,H}$ is 
\beq \label{BSDEYH}
Y_t^{0,H} &=& H_T^0 + \frac{1}{p} \ln \gamma_T^0 + \int_t^T  f^{0,H}(r,Y_r^{0,H},Z_r^{0,H}) dr  - \int_t^T Z_r^{0,H} dW_r,
\enq
where the generator $f^{0,H}$ is the  $\Oc(\F)\otimes\Bc(\R)\otimes\Bc(\R)$-measurable map defined by
\beq
f^{0,H}(t,y,z) &=&   -   \frac{b_t^0}{\sigma_t^0}  z  -  \frac{1}{2p}   \Big(\frac{b_t^0}{\sigma_t^0} \Big)^2 \label{genf0H} \\
& & \;\; +  \frac{p}{2} \inf_{a\in A_0} 
 \Big| \Big( z+ \frac{1}{p} \frac{b_t^0}{\sigma_t^0} \Big) -  a \sigma_t^0  + \frac{2}{p}  U(y) \int_E U(ae-Y_t^{1,H}(t,e) \eta(de)  \Big|^2.  \nonumber
\enq

The solution to the optimal investment problem without defaultable claim is obtained similarly as  for the case with claim, by  considering 
$H$ $=$ $0$.  We thus have $V_0(x)$ $=$ $R(x-Y_0^0)$, where the BSDE associated to $Y^0$ is given by
\beqs
Y_t^{0} &=&  \frac{1}{p} \ln \gamma_T^0 + \int_t^T  f^0(r,Y_r^{0},Z_r^{0}) dr  - \int_t^T Z_r^{0} dW_r,
\enqs
with a generator $f^0$ as in \reff{genf0H} for $H$ $=$ $0$, i.e. $Y^{1,H}$ replaced by $Y^1$ solution to the BSDE 
\beqs
Y_t^{1}(\theta,e) &=&   \frac{1}{p} \ln \gamma_T(\theta,e)  + \int_t^T f^{1}(r,Z_r^{1},\theta,e) dr - \int_t^T Z_r^{1} dW_r.
 \enqs


Finally, the indifference price is given by 
\beqs
\pi &=& Y_0^{0,H} - Y_0^0. 
\enqs

\begin{Remark}
{\rm   Notice that the quadratic generator $f^{0,H}$ in \reff{genf0H} of the BSDE \reff{BSDEYH} is not standard due to the additional 
term arising from the integral gain involving $Y^{1,H}$. However, one can prove existence and uniqueness of this BSDE and obtain a 
verification theorem relating the solution of this BSDE to the original value function by choosing a suitable set of admissible controls $\Ac_\G$ $=$ 
$\Ac_\F^0\times\Ac_\F^1$.  The details are provided in the companion paper \cite{JKP}. 
Actually, in this related paper, we consider a multi-dimensional extension of the above model with assets subject to successive  counterparty 
default times, and we apply the $\F$-decomposition method for solving  the indifference pricing of defaultable claims, including  credit derivatives such as $k$-th default swap.  
}
\end{Remark}

\subsection{Optimal investment under bilateral counterparty risk}

We consider a portfolio with two names, each one  subject to an external counterparty default, but also to the default of the other one due to 
a contagion effect.    
We denote by $S^1$ and $S^2$ the value process of these two names, by $\tau_1$ and $\tau_2$ their default  times, not necessarily ordered, 
and by $\hat\tau_1$ $=$ $\min(\tau_1,\tau_2)$, $\hat\tau_2$ $=$ $\max(\tau_1,\tau_2)$. 
Once the name $i$ defaults at random time $\tau_i$, meaning that the value of $S^i$ drops to zero, it also incurs a jump (drop or gain) on the other value process  $S^j$, $i,j$ $\in$ $\{1,2\}$, $i$ $\neq$ $j$.  

The reference filtration $\F$ is the filtration generated by a two-dimensional Brownian motion $W$ $=$ $(W^1,W^2)$, driving the evolution of  the names in absence of defaults,  
and the global market information is represented by  $\G$ $=$ $\F$ $\vee$ $\D^1$ $\vee$ $\D^2$, with $\D^i$ $=$ $(\Dc_t^i)_{t\geq 0}$, 
$\Dc_t^i$ $=$ $\cap_{\eps>0} \sigma(1_{\tau_i\leq s}, s \leq t+\eps)$, $i$ $=$ $1,2$.

The $\G$-adapted value processes  $S^i$ of names $i$ $=$ $1,2$,  are  given by 
\beqs
S_t^i &=& S_t^{i,0} 1_{t<\hat\tau_1} + S_t^{i,j}(\tau_j) 1_{\tau_j\leq t < \tau_i}, \;\;\; t \geq 0, \;\; i,j =1,2, \; i \neq j, 
\enqs
where $S^0$ $=$ $(S^{1,0},S^{2,0})$ is the vector price process of the two names in absence of any default,  governed by 
\beqs
dS_t^{0} &=& {\rm diag}(S_t^{0}) \big( b_t^{0} dt + \sigma_t^{0} dW_t\big),
\enqs
$b^0$ $=$ $(b^{1,0},b^{2,0})$ is $\F$-adapted, $\sigma^0$ is the $2\times 2$-diagonal $\F$-adapted matrix with diagonal diffusion coefficients 
$\sigma^{1,0}$ $>$ $0$, $\sigma^{2,0}$ $>$ $0$, and the indexed process $S^{i,j}(\theta_j)$, $\theta_j$ $\in$ $\R_+$, representing the 
price process of name $i$ after the default of name $j$ at time $\theta_j$, is given by 
\beqs
dS_t^{i,j}(\theta_j) &=& dS_t^{i,j}(\theta_j) \big( b_t^{i,j}(\theta_j) dt + \sigma_t^{i,j}(\theta_j) dW_t^i \big), \;\; t \geq \theta_j, \\
S_{\theta_j}^{i,j}(\theta_j) &=& S_{\theta_j}^{i,0}.(1 + e^{i,j}),
\enqs
where $e^{i,j}$ represents the proportional jump induced by the default of name $j$ on name $i$, and assumed  constant for simplicity and 
valued in $(-1,\infty)$ . The coefficients $b^{i,0}$, $\sigma^{i,0}$ $>$ $0$ 
are $\F$-adapted processes, and $b^{i,j}$, $\sigma^{i,j}$ $>$ $0$ are in $\Oc^1_{\F}(\R_+)$.  

The trading strategy of the investor is  a  $\G$-predictable measurable process $\alpha$ representing the fraction of wealth 
invested in the two names. 
It is then decomposed in four components: the first component $\alpha^0$  is a pair  of  $\F$-predictable processes  representing the fraction invested in the two names before any default, the second component $\alpha^{1,1}$ is an indexed $\F$-predictable process representing the fraction invested in the name $2$ when the name $1$ defaults, the third component $\alpha^{1,2}$ 
is an indexed  $\F$-predictable process representing the fraction invested in the name $1$ when the name $2$ defaults, and the fourth component is zero when both names default. The wealth process of the investor is then given by 
\beqs
X_t &=& X_t^0 1_{t<\hat\tau_1} + X_t^{1,1}(\tau_1) 1_{\tau_1\leq t<\tau_2} + X_t^{1,2}(\tau_2) 1_{\tau_2\leq t<\tau_1} + 
X_t^2(\tau_1,\tau_2) 1_{t\geq\hat\tau_2}, 
\enqs
where $X^0$ is the wealth process before any default, governed by 
\beqs
dX_t^0 &=& X_t^0 (\alpha_t^0)'{\rm diag}(S_t^0)^{-1}dS_t^0 \\
&=&   X_t^0 \big(\alpha_t^{0}.b_t^{0} dt  +  (\alpha_t^{0})'\sigma_t^{0} dW_t \big),
\enqs
$X^{1,1}(\theta_1)$ is the wealth indexed process after default of name $1$,  governed by
\beqs
dX_t^{1,1}(\theta_1) &=& X_t^{1,1}(\theta_1) \alpha_t^{1,1}(\theta_1) \frac{dS_t^{2,1}(\theta_1)}{S_t^{2,1}(\theta_1)}, \; t \geq \theta_1 \\
X_{\theta_1}^{1,1}(\theta_1) &=& X_{\theta_1}^0. ( 1 + \alpha_{\theta_1}^{0}.(-1,e^{2,1})),
\enqs
$X^{1,2}(\theta_2)$ is the wealth indexed process after default of name $2$, governed by 
\beqs
dX_t^{1,2}(\theta_2) &=& X_t^{1,2}(\theta_2) \alpha_t^{1,2}(\theta_2) \frac{dS_t^{1,2}(\theta_2)}{S_t^{1,2}(\theta_2)}, \; t \geq \theta_2 \\
X_{\theta_2}^{1,2}(\theta_2) &=& X_{\theta_2}^0 . \big( 1 +  \alpha_{\theta_2}^{0}.(e^{1,2},-1)\big),
\enqs
and $X^2(\theta_1,\theta_2)$ is the wealth indexed process after both defaults,  hence constant after $\theta_1\vee\theta_2$, and  then given by
\beqs
X_t^2(\theta_1,\theta_2) &=& \left\{
					       \begin{array}{cc}	
					       X_{\theta_2}^{1,1}(\theta_1) . \big(1 - \alpha_{\theta_2}^{1,1}(\theta_1)\big),& \;  \theta_1 \leq \theta_2 \leq t \\
					       X_{\theta_1}^{1,2}(\theta_2) . \big(1 - \alpha_{\theta_1}^{1,2}(\theta_2)\big),& \;  \theta_2 < \theta_1 \leq t
					       \end{array}
					       \right. 	
\enqs
In order to ensure that the wealth process is strictly positive, we assume that  $\alpha^0$ is valued in a closed subset $A_0$ $\subset$ 
$\{ a \in \R^2:  1+ a.(-1,e^{2,1}) > 0, \;\mbox{ and }  1+ a.(e^{1,2},-1) >0 \}$, and $\alpha^{1,1}$, $\alpha^{1,2}$ are valued respectively in 
closed subsets $A_{1,1}$, $A_{1,2}$ $\subset$ $(-\infty,1)$.

We are next given a utility function $U$ on $\R_+$, over a finite horizon $T$, and we consider the optimal investment problem
\beq \label{defV0bi}
V_0(x) &=& \sup_{\alpha\in\Ac_{\G}} \E\big[ U(X_T^{x,\alpha})\big]. 
\enq

We use the $\F$-decomposition method of Section \ref{secmultiple} for the resolution of \reff{defV0bi}. From Theorem \ref{thm2}, the value function 
$V_0$ is  obtained via the following backward induction formula:
\beqs
V_2(x,\theta_1,\theta_2) &=& U(x) \E\big[ \gamma_T(\theta_1,\theta_2) \big| \Fc_{\theta_1\vee\theta_2} \big] 
\; := \; U(x) \bar\gamma(\theta_1,\theta_2) \\
V_{1,1}(x,\theta_1) &=&  \esssup_{\alpha^{1,1}\in\Ac^{1,1}_{\F}} \E \Big[ U(X_T^{1,1,x}) \gamma_T^{1,1}(\theta_1) + 
\int_{\theta_1}^T U\big(X_{\theta_2}^{1,1}.\big(1 - \alpha_{\theta_2}^{1,1}\big)\big)  \bar\gamma(\theta_1,\theta_2) 
d \theta_2 \Big | \Fc_{\theta_1} \Big] \\
V_{1,2}(x,\theta_2) &=&  \esssup_{\alpha^{1,2}\in\Ac^{1,2}_{\F}} \E \Big[ U(X_T^{1,2,x}) \gamma_T^{1,2}(\theta_2) + 
\int_{\theta_2}^T U\big(X_{\theta_1}^{1,2}.\big(1 - \alpha_{\theta_1}^{1,2}\big)\big)  \bar\gamma(\theta_1,\theta_2) 
d \theta_1 \Big | \Fc_{\theta_1} \Big] \\
V_0(x) &=& \sup_{\alpha^0\in\Ac_{\F}^0} \E\Big[ U(X_T^{0,x})\gamma_T^0   \\
& &\;\; + \;  \int_0^T   V_{1,1}\big( X_{\theta}^0.\big( 1 + \alpha_{\theta}^{0}.(-1,e^{2,1})\big),\theta) 
+ V_{1,2}\big( X_{\theta}^0.\big( 1 +  \alpha_{\theta}^{0}.(e^{1,2},-1)\big),\theta\big) \; d\theta \Big]. 
\enqs

In the sequel, we  consider power utility functions
\beqs
U(x) &=& \frac{1}{p}  x^p, \;\;\; x \geq 0, \; p < 1, \; p \neq 0, 
\enqs
and we use dynamic programming and BSDE methods in the $\F$-filtration to solve the above stochastic control problems.  
The value functions $V_{1,1}$ and $V_{1,2}$ are then in the form
\beqs
V_{1,1}(x,\theta_1) \; = \; U(x) Y^{1,1}_{\theta_1}(\theta_1), & &  V_{1,2}(x,\theta_2) \; = \; U(x) Y^{1,2}_{\theta_2}(\theta_2),
\enqs
where $Y^{1,1}(\theta_1)$ and $Y^{1,2}(\theta_2)$ are solutions to the BSDEs: 
\beqs
Y_t^{1,1}(\theta_1) &=& \gamma_T^{1,1}(\theta_1)  + \int_t^T f_{1,1}(r,Y_r^{1,1}(\theta_1),Z_r^{1,1}(\theta_1),\theta_1) dr 
- \int_t^T Z_r^{1,1}(\theta_1) dW^2_r, \\
Y_t^{1,2}(\theta_2) &=& \gamma_T^{1,2}(\theta_2)  + \int_t^T f_{1,2}(r,Y_r^{1,2}(\theta_2),Z_r^{1,2}(\theta_2),\theta_2) dr, 
- \int_t^T Z_r^{1,2}(\theta_2) dW^1_r 
\enqs
with  generators
\beqs
f_{1,1}(t,y,z,\theta_1) &=& p  \; \sup_{a \in A_{1,1}} \Big[  \big(b_t^{2,1}(\theta_1)y + \sigma_t^{2,1}(\theta_1)z\big) a - 
\frac{1-p}{2} y  |  \sigma_t^{2,1}(\theta_1) a |^2 \\
& & \hspace{25mm} + \;  \bar\gamma(\theta_1,t) \frac{(1-a)^p}{p} \Big] \\
f_{1,2}(t,y,z,\theta_2) &=& p \;  \sup_{a \in A_{1,2}} \Big[  \big(b_t^{1,2}(\theta_2)y + \sigma_t^{1,2}(\theta_2)z\big) a - 
\frac{1-p}{2} y  |  \sigma_t^{1,2}(\theta_2) a |^2 \\
& & \hspace{25mm} + \;  \bar\gamma(t,\theta_2) \frac{(1-a)^p}{p} \Big]. 
\enqs
Finally, we have 
\beqs
V_0(x) &=& U(x) Y_0,
\enqs
where $Y^0$ is the solution to the BSDE
\beqs
Y_t^0 &=& \gamma_T^0 +  \int_t^T f_0(r,Y_r^0,Z_r^0) dr - \int_t^T Z_r^0. dW_r,
\enqs
with a generator
\beqs
f_0(t,y,z) &=& p \;  \sup_{a \in A_0} \Big[  \big( y b_t^0 + \sigma_t^0 z).a - \frac{1-p}{2} y |\sigma_t^0 a|^2 \\
& & \hspace{13mm} + \; Y_t^{1,1}(t) \frac{(1+a.(-1,e^{2,1}))^p}{p} +  Y_t^{1,2}(t) \frac{(1+a.(e^{1,2},-1))^p}{p} \Big]. 
\enqs
The details and rigorous mathematical treatment of the above derivation are studied in \cite{JKP}, where we prove 
the existence and uniqueness of the  solutions to these BSDEs, and that they are  
indeed related to the original value functions of our optimal investment problem.

\end{document}